\documentclass[preprint,12pt]{elsarticle}

\usepackage{amssymb}
\usepackage{amsmath}
 \usepackage{amsthm}

\journal{Computer Aided Geometric Design}
\usepackage{txfonts}
\usepackage{amsmath}
\usepackage{amssymb}
\usepackage{graphicx}
\usepackage{hyperref}
\usepackage{booktabs}
\usepackage{multirow}
\hypersetup{
    colorlinks=true,
    linkcolor=blue,
    filecolor=magenta,      
    urlcolor=cyan,
    pdftitle={Overleaf Example},
    pdfpagemode=FullScreen,
    }
\newtheorem{theorem}{Theorem}[section]
\newtheorem{lemma}[theorem]{Lemma}
\newtheorem{proposition}[theorem]{Proposition}

\theoremstyle{definition}
\newtheorem{definition}[theorem]{Definition}
\newtheorem{example}[theorem]{Example}

\theoremstyle{remark}
\newtheorem{remark}[theorem]{Remark}

\numberwithin{equation}{section}

\newcommand{\placket}[1]{\left({#1}\right)}
\newcommand{\vlacket}[1]{\left\{{#1}\right\}}
\newcommand{\blacket}[1]{\left[{#1}\right]}

\newcommand{\imaginary}[0]{\sqrt{-1}\,}

\newcommand{\note}[1]{\texttt{ <<Note:{#1}>> }}

\begin{document}

\begin{frontmatter}



\title{The Klein-Lie W-curves as a new class of aesthetic curves based on self-affinity} %% Article title

 \author[label1]{Shun Kumagai}
 \affiliation[label1]{organization={Faculty of Engineering, Hachinohe Institute of Technology},
             addressline={88-1 Myo Ohbiraki, Hachinohe},
            city={Aomori},
            postcode={031-8501},
            country={Japan}}
 \author[label2]{Kenji Kajiwara}
 \affiliation[label2]{organization={Institute of Mathematics for Industry, Kyushu University},
            addressline={774 Motooka, Nishi-ku},
            city={Fukuoka},
            postcode={819-0395},
            country={Japan}}

%% Abstract

\begin{abstract}
In this paper, we consider planar curves and present a family of planar curves characterized by a symmetry called the extendable self-affinity (ESA). The ESA has been recognized through the investigation of the symmetry of the log-aesthetic curve (LAC), which has been studied as a reference for designing aesthetic shapes in CAGD and regarded as an analog of Euler's elastica in similarity geometry. 
We investigate the ESA and show that it gives rise to affine geometry and the Klein-Lie W-curve. 
An asymptotic behavior of the logarithmic curvature graph of a subclass of the W-curve is proven to be linear as an invariant property under affine transformations, by evaluating similarity curvature. 
Based on the observations on the self-affinity and the logarithmic curvature graph, we propose the Klein-Lie W-curve as a new family of ``aesthetic curves" in affine geometry, to be considered together with the LAC. 

\end{abstract}

\if0 %version 20250428
The log-aesthetic curve is a family studied as a reference for designing aesthetic shapes in CAD systems. 
It is desirable for CAGD to give a 
unified description of ``aesthetic shapes" in the framework of Klein geometry.

In this paper, we consider curves in equiaffine geometry. 
Then, we present a new family of plane curves characterized by a self-affinity called the extendable self-affinity, which has been 
recognized through the investigation of the symmetry of the log-aesthetic curve.

Our new family includes the logarithmic spiral, a special case of the log-aesthetic curves, and the quadratic curves, which implies that the family can be regarded as an alternate family of "aesthetic curves".
We discuss the relationship with the known two self-affinities, namely, the Miura self-affinity and the Harada self-affinity for a unified description of "aesthetic curves".
\fi



%% Keywords

\end{frontmatter}




\section{Introduction}
%\subsection{Background and overview}
The log-aesthetic curve (LAC) is a family of planar curves studied as a reference for designing aesthetic shapes in CAGD \cite{Harada1999, Kanaya2003en, Yoshida2006, Miura2007, Yoshida2012, Gobi20143, Graiff2023}. In this paper, we consider an extension of the LAC based on the symmetry that we call \textit{self-affinity} \cite{Harada1995, Miura2006, KK}. 

The LAC has been identified and proposed by Harada et al. \cite{Harada1995} through the effort of extracting the common properties of the planar curves which car designers regard as \textit{aesthetic}. 
They considered the arc length of curve segments in a specific range of curvature radius, plotted the log-log histogram of the former versus the latter, which is called the \textit{logarithmic curvature histogram (LCH, also
known as the logarithmic distribution diagram of curvature, LDDC)}. 
They pointed out that the ``aesthetic" curves are identified by the linear tendencies of their LCHs, and called them \textit{monotonic rhythm curves}. 
Later, Miura \cite{Miura2006} considered the continuum limit of the LCH called the \textit{logarithmic curvature graph (LCG)} to arrive at the definition of the LAC from the exact linearity of the LCG. It is given by the formula
\begin{align}\rho(s)=
    \begin{cases}
        (\xi s +\eta)^{\frac{1}{\alpha}} &(\alpha\neq 0),
        \\
       e^{\xi s +\eta} &(\alpha= 0),
    \end{cases}\quad\xi,\eta\in\mathbb{R}, 
\end{align}
where $\rho$ is the (Euclidean) curvature radius, $s$ is the arc length and $\alpha\in\mathbb{R}$ is the slope of the linear LCG. 
%It is a linear reparametrization of  Savelov's pseudo-spiral \cite{Savelov, Ziatdinov}. %Savelov, 1960 / Gray 1997

The LAC has been discussed in the framework of similarity geometry in \cite{Inoguchi2018, InoguchiMiura2019, Inoguchi2023}, where it is the shape invariant curve with respect to integrable deformation in similarity geometry. 
It is also shown that the LAC admits the variational formulation in terms of the fairing energy functional. 
Those characterizations suggest that the LAC is regarded as an analogue of Euler's elastica \cite{Elastica} in similarity geometry. 

On the other hand, the first consideration %foundation 
\cite{Harada1995}  of monotonic rhythm curves by Harada et al.\ entails a concept of symmetry of curves called the \textit{self-affinity}. 
Miura \cite{Miura2006} modified the definition of self-affinity to fit the LAC.  
It has been proved in \cite{KK} that the Miura self-affinity (MSA) characterizes the LAC. 
It has also been shown that the Harada self-affinity (HSA) does not characterize the LAC but actually \textit{parabolas}, which is the zero-curvature curves in equiaffine geometry.
The HSA can be relaxed into the \textit{extendable self-affinity (ESA)} \cite{KK} to capture a wider family of curves with reference to the MSA. 
Then, it is shown that the ESA gives rise to quadratic curves, which are the constant-curvature curves in the equiaffine geometry and are also fundamental components in CAD systems \cite{CAGD, Bezier}.

Motivated by the above results, we consider the ESA of quadratic curves in a weaker sense, which gives rise to the ESA in affine geometry. 
It results in a subfamily of so-called the \textit{Klein-Lie W-curve} \cite{LieKlein, Blaschke}, which we call \textit{affine W-curve}. 
It is the family of curves with constant affine curvature and includes the graphs of power function, exponential function, $y=x\log x$, quadratic curves and \textit{logarithmic spiral}. 
Noticing that logarithmic spiral also belongs to the LAC, this suggests that the affine W-curve and the LAC should be considered in a unified manner as a new family of aesthetic curves.

Notably, the affine W-curve asymptotically behaves like the LAC in terms of the linearity of the LCG. 
We present several examples showing that the LCGs of affine W-curves exhibit asymptotic linearity near their asymptotic ends, which suggests that the asymptotic slope is preserved under affine deformations.
We subsequently prove that this phenomenon holds universally for any non-degenerate affine W-curve admitting invariant asymptotic directions under the ESA. 
Specifically, we establish this result by evaluating the asymptotic expansion of the similarity curvature and analyzing the action of the affine transformation group on Euclidean and similarity geometric quantities.
The resulting asymptotic behavior also strongly motivates a unified treatment of the affine W-curve and the LAC from both algebraic and analytic viewpoints.

%

%It is pointed out \cite{Harada1999} that the concept of self-affinity is an important for artisticity of curves from the viewpoint of design;  
%It is claimed in \cite{Harada1999} that the designers aimed to draw the curve as having the self-affinity. 
%A classification of curves in visual language into five types in terms of 

Throughout this paper, we consider sufficiently smooth, parametric planar {curves}. %with $C^2$-continuity. 
%Let $C^2(I,{\mathbb{R}^2})$ denote the set of such curves with domain $I\subset \mathbb{R}$. 
For a curve $\gamma(t):\mathbb{R}\rightarrow{\mathbb{C}}\cong \mathbb{R}^2$, a \textit{reperametrization} of $\gamma(t)$ is a smooth, one-to-one map $w=w(t):\mathbb{R}\rightarrow \mathbb{R}$.
In this case, we say that $\gamma$ admits a parametrization $\gamma(w)$ in the parameter $w=w(t)$, where $\gamma(w):=(\gamma\circ w^{-1})(w)$.
For a given parametrization, we fix a base parameter represented with a subscript zero, such as $t_0$. 
Quantities associated with a curve $\gamma$ may be explicitly indicated by the subscript $\gamma$. 
Unless otherwise noted, subscripts denote differentiation with respect to the corresponding variable.

\section{Klein geometries, log-aesthetic curves and self-affinities}

In this section, we first give a brief review of the basic framework of planar curves in Klein geometries \cite{Inoguchi2023, Equiaffine} that is used in this paper. 
We then introduce the LAC and self-affinities as its characterizing properties.    
We give a summary of discussions about the HSA, the MSA and the ESA, in particular, how a new family of curves arises in affine geometry.

\subsection{Euclidean and similarity geometries}
If $\gamma_t(t)\neq0$ for any $t$, we have the \textit{arc length} parameter
$s=s(t):=\int_{t_0}^t|\gamma_t(t)|dt$
for which $|\gamma_s|=1$ and $s_0=s(t_0)=0$ hold. 
The fundamental theorem of curves states that 
for a given non-negative, smooth function $\kappa^\mathrm{E}(s):\mathbb{R}\rightarrow \mathbb{R}$, 
the \textit{Frenet formula} 
\begin{align}
        \Phi^\mathrm{E}_s=\Phi^\mathrm{E}\begin{pmatrix}
            0&-\kappa^\mathrm{E}\\ \kappa^\mathrm{E}&0
        \end{pmatrix},\quad \Phi^\mathrm{E}=(\gamma_s,\imaginary\gamma_{s}),\label{eqFrenet}
\end{align}
    has a unique solution $\gamma(s)$ up to the Euclidean motion group 
    \begin{align}
        G^\mathrm{E}:=\{F:\mathbb{C}\to \mathbb{C}\mid F(z)= Az+b,\ A\in \mathrm{SO}(2),\ b\in{\mathbb{R}^2}\},
    \end{align} where a matrix $A$ acts linearly on $\mathbb{C}\cong \mathbb{R}^2$.
    The function $\kappa^\mathrm{E}$ is called the \textit{(Euclidean) curvature} of $\gamma$. 
    
In the same setting 
$\gamma_t\neq0$, we also have the \textit{similarity arc length (turning angle)} parameter 
$
    \theta=\theta(s):=\theta_0+\int_0^s \kappa^\mathrm{E}(s)ds. 
$
Then, %the following similarity analog of the fundamental theorem of curves \cite{sim} states that 
%Fujioka-Inoguchi 2007 Progress in Math
it follows that for a given non-negative, smooth function $\kappa^\mathrm{sim}(\theta):\mathbb{R}\rightarrow \mathbb{R}$  the \textit{similarity Frenet formula}
\begin{align}\label{simFrenet}
        \Phi^\mathrm{sim}_\theta=\Phi^\mathrm{sim}\begin{pmatrix}
            -\kappa^\mathrm{sim}&-1\\ 1&-\kappa^\mathrm{sim}
        \end{pmatrix},\quad \Phi^\mathrm{sim}=(\gamma_\theta,\imaginary\gamma_{\theta}),
\end{align}
has a unique solution $\gamma(\theta)$ up to the similarity transformation group 
\begin{align}
    G^\mathrm{sim}:=\{F:\mathbb{C}\to \mathbb{C}\mid F(z)= Az+b,\ A\in \mathbb{R}^+\times \mathrm{SO}(2),\ b\in{\mathbb{R}^2}\}.
\end{align}
The function $\kappa^\mathrm{sim}={\kappa^\mathrm{E}_\theta}/{\kappa^\mathrm{E}}={\kappa^\mathrm{E}_s}/{(\kappa^\mathrm{E})}^2$ is called the \textit{similarity curvature} of $\gamma$.

\subsection{Equiaffine geometry}
A curve $\gamma(t)$ is called \textit{non-degenerate} if $\det(\gamma_t(t),\gamma_{tt}(t))\neq0$. For a non-degenerate curve $\gamma(t)$, we have the \textit{equiaffine arc length} parameter $
    u=u(t):=u_0+\int_{t_0}^t\det(\gamma_t(t),\gamma_{tt}(t))^\frac{1}{3}dt,$
for which $\det(\gamma_u,\gamma_{uu})=1$ follows. 
For a given smooth function $\kappa^\mathrm{SA}(u):\mathbb{R}\rightarrow \mathbb{R}$, 
the \textit{equiaffine Frenet formula} 
\begin{align}\label{SAFrenet}
        \Phi^\mathrm{SA}_u=\Phi^\mathrm{SA}\begin{pmatrix}0&-\kappa^\mathrm{SA}\\ 1&0
        \end{pmatrix},\quad \Phi^\mathrm{SA}=(\gamma_u,\gamma_{uu}),
\end{align}
has a unique solution $\gamma(u)$ up to the equiaffine transformation group 
\begin{align}
    G^\mathrm{SA}:=\{F:\mathbb{C}\to \mathbb{C}\mid F(z)= Az+b,\ A\in \mathrm{SL}(2,\mathbb{R}),\ b\in{\mathbb{R}^2}\}.
\end{align}
The function $\kappa^\mathrm{SA}$ is called the \textit{equiaffine curvature} of $\gamma$.
\begin{proposition}\label{kappaSA}
    The equiaffine curvature is represented by 
\begin{align}\label{kSASchwarz}
    \kappa^\mathrm{SA}&=-\frac{\gamma_{uuu}}{\gamma_u}=-\frac{1}{2}S_u(\tan \theta)=-\frac{1}{2}\placket{\frac{(\tan\theta )_{uu}}{(\tan\theta )_u}}_u+\frac{1}{4}\placket{\frac{(\tan\theta )_{uu}}{(\tan\theta )_u}}^2,
\end{align}
where $\tan \theta=\mathrm{Im}(\gamma_t)/\mathrm{Re}(\gamma_t)$ is the slope of the tangent vector, and $S_u$ denotes the Schwarzian derivative with respect to $u$. 
\end{proposition}
\begin{proof}
    Denote $\gamma(u)=x(u)+\imaginary y(u)$ and $z(u)=\tan\theta =y_u/x_u$. 
    From equiaffine Frenet formula \eqref{SAFrenet}, we have $\gamma_{uuu}=-\kappa^\mathrm{SA}\gamma_u$, so that $x_{uuu}=-\kappa^\mathrm{SA}x_u$. 
    Then, we show $x_{uuu}=-S_u(z) x_u/2$, from which we identify $\kappa^\mathrm{SA}(u)=S_u(z)/2$ by the uniqueness of solution of equiaffine Frenet formula \eqref{SAFrenet}.
    We have
    \begin{align}
        %z_u&=\frac{y}{x}\placket{\frac{y_u}{y}-\frac{x_u}{x}}, \\
         %z_{uu}&=\frac{2x_uy}{x^2}\placket{\frac{y_u}{y}-\frac{x_u}{x}},\\
        \frac{z_{uu}}{z_u}&=\frac{\dfrac{2x_{uu}y_u}{x_u^2}\placket{\dfrac{y_{uu}}{y_u}-\dfrac{x_{uu}}{x_u}}}{\dfrac{y_u}{x_u}\placket{\dfrac{y_{uu}}{y_u}-\dfrac{x_{uu}}{x_u}}}=\frac{2x_{uu}}{x_u}. \label{eq1}
    \end{align}
    Differentiating \eqref{eq1} by $u$ yields 
    \begin{align}
        \placket{\frac{z_{uu}}{z_u}}_u&=-\frac{2x_{uuu}}{x_u}+\frac{2x_{uu}^2}{x_u^2}=-\frac{2x_{uuu}}{x_u}+\frac{1}{2}\placket{\frac{z_{uu}}{z_u}}^2,
        \end{align}
        from which we have
        \begin{align}
        \frac{x_{uuu}}{x_u}&=-\frac{1}{2}\vlacket{\placket{\frac{z_{uu}}{z_u}}_u-\frac{1}{2}\placket{\frac{z_{uu}}{z_u}}^2}=-\frac{1}{2}S_u(z).\label{eq2}
    \end{align}
    %We have conclusion from \eqref{eq2}
\end{proof}
For  $F=(z\mapsto Az+b)\in G^\mathrm{SA}$, the Jacobian $DF$ is $A\in \mathrm{SL}(2,\mathbb{R})$. 
The action of $F$ preserves the equiaffine arc length $u(t)$ of a non-degenerate curve $\gamma(t)$.
The action of $F$ on the slope of tangent vector is governed by the M\"obius transformation 
\begin{align}\label{Mobslope}
    \tan\theta_{F\gamma}=\frac{\mathrm{Im}(F\gamma_t)}{\mathrm{Re}(F\gamma_t)}=\frac{cx_t+dy_t}{a x_t+by_t}=\frac{c+d\tan\theta_\gamma}{a+b\tan\theta_\gamma},\quad DF=A=\begin{pmatrix}
        a&b\\c&d
    \end{pmatrix},
\end{align}
where $\theta_\gamma$ and $\theta_{F\gamma}$ denote the turning angles of $\gamma$ and $F\gamma$, respectively. 
By the well-known fact that the Schwarzian derivative is invariant under M\"obius transformations, Proposition \ref{kappaSA} implies that $\kappa^\mathrm{SA}(u)$ is invariant under the action of $F \in G^\mathrm{SA}$. 
%Therefore, the Schwarzian derivative appears naturally for describing the invariant $\kappa^\mathrm{SA}$ under equiaffine transformations. 
As a trivial case, the following results \cite{Blaschke} are known. 
\if0
It is represented by the formula \cite{Equiaffine}
    \begin{align}
        \kappa^\mathrm{SA}&=-\frac{\gamma_{uuu}}{\gamma_u}=(\kappa^\mathrm{E})^{\frac{4}{3}}+\frac{1}{3}(\kappa^\mathrm{E})^{-\frac{5}{3}}\kappa^\mathrm{E}_{ss}-\frac{5}{9}(\kappa^\mathrm{E})^{-\frac{8}{3}}(\kappa^\mathrm{E}_s)^2.
    \end{align}
\fi
\begin{lemma}\label{quad}
    %A curve in equiaffine geometry (the solution of equiaffine Frenet formula \eqref{SAFrenet}) 
    A non-degenerate curve has constant equiaffine curvature if and only if it is a quadratic curve.
    More precisely, it is either
    \begin{enumerate}
        \item a parabola ($\kappa^\mathrm{SA}=0$), 
        \item an ellipse ($\kappa^\mathrm{SA}>0$), 
        \item a hyperbola ($\kappa^\mathrm{SA}<0$). 
    \end{enumerate}
\end{lemma}

\if0
\begin{proof}
    %By solving the differential equation $\gamma_{uuu}=-\kappa^\mathrm{SA}\gamma_u$ one obtains %the conclusion. 
%\end{proof}
We suppose the following presentation of quadratic curves:
\begin{enumerate}
    \item $\gamma(t)=t+\imaginary t^2$ (a parabola),
    \item $\gamma(t)=\cos t +\imaginary\sin t$ (an ellipse, equivalently a circle),
    \item $\gamma(t)=\cosh t +\imaginary \sinh t$ (a hyperbola).
\end{enumerate}
One can calculate each of the above yields a constant equiaffine curvature $\kappa^\mathrm{SA}$, which is invariant under the action of $G^\mathrm{SA}$. 
Next, for $a>0$ and a curve $\gamma(u)$ with constant curvature $\kappa^\mathrm{SA}(u)$, the equiaffine length $u^a(u)$ of the scaled curve $\gamma_a(u):=a\gamma(u)$ is given by 
\begin{align}
    u^a(u)=u^a_0+\int_{u_0} ^u \det (a\gamma_u, a\gamma_{uu})du& =u^a_0+a^2\int_{u_0} ^u \det (\gamma_u, \gamma_{uu})du
    \\&=u^a_0+ a^2(u-u_0). 
\end{align}
Thus the equiaffine curvature is calculated by $\kappa^\mathrm{SA}_{a\gamma}=-\frac{a\gamma_{u^au^au^a}}{a\gamma_{u^a}}=\frac{1}{a^4}\kappa^\mathrm{SA}_\gamma$ and any constant $\kappa^\mathrm{SA}$ is represented by quadratic curves. 
Finally, by the uniqueness of solutions of linear differential equation \eqref{SAFrenet}, we conclude the claim. 
\end{proof}
\fi

For a general equiaffine curvature $\kappa^\mathrm{SA}(u)$, the fundamental theory of linear ordinary differential equations applied to the equiaffine Frenet formula \eqref{SAFrenet} yields the following representation formula.
\begin{lemma}\label{lemma:repformula}
    Let $z(u)=f(u),g(u)$ be the fundamental solutions of the differential equation $z_{uu}(u)+\kappa^\mathrm{SA}(u)z(u)=0$ with $W(f,g):=\det\begin{pmatrix}
        f&f_u\\g&g_u
    \end{pmatrix}=1$. 
    Then, the solution $\gamma(u)$ of the equation $\gamma_{uuu}=-\kappa^\mathrm{SA}(u)\,\gamma_u$ is up to $G^\mathrm{SA}$ given by 
    \begin{align}\label{repformula}
        \gamma(u)-\gamma(u_0)=\int^u_{u_0}\vlacket{f(u)+\imaginary g(u)}du.
    \end{align}    
\end{lemma}

\subsection{Affine gemoetry}
Under the setting $\det(\gamma_{u},\gamma_{uu})=1$,
we consider invariants under the action of the \textit{affine group} 
\begin{align}\label{defGA}
    G^\mathrm{A}:=\vlacket{F:\mathbb{C}\to \mathbb{C}\mid F(z)= Az+b,\ A\in \mathrm{GL}^+(2,\mathbb{R}), b\in \mathbb{R}^2}. 
\end{align}
We may define the \textit{affine arc length} parameter $\sigma(u)$ and \textit{affine curvature} $\chi(u)$ \cite[p.58]{Encyclopedia} as
\begin{align}\label{sigma}
    \sigma(u)&:=\int {\kappa^\mathrm{SA}(u)}^{\frac{1}{2}}du,
    \\
    \chi(u)&:={\kappa^\mathrm{SA}(u)}^{-\frac{3}{2}}\kappa^\mathrm{SA}_u(u),
\end{align}
respectively.
    We note that the M\"obius transformation formula \eqref{Mobslope} remains valid for any $F \in G^\mathrm{A}$ as well as orientation-reversing affine transformations.  
    However, under the action of $F \in G^\mathrm{A}$, the equiaffine quantities $u, \kappa^\mathrm{SA}$ are no longer invariant, and the affine quantities $\sigma,\chi$ play the role of affine invariants as follows. 
\begin{lemma}\label{affinv}
    Let $\gamma$ be a non-degenerate curve and $F\in G^\mathrm{A}$. Denote the affine arc lengths of $\gamma, F\gamma$ by $\sigma_\gamma,\sigma_{F\gamma}$ and the affine curvatures of $\gamma, F\gamma$ by $\chi_\gamma,\chi_{F\gamma}$, respectively. Then, the following holds. 
    %Under the above setting, the following holds.
    \begin{enumerate}
        \item The affine curvature $\chi_\gamma$ takes value in $\mathbb{R}\cup \imaginary\mathbb{R}$. Furthermore, $\mathrm{sign}(\chi^2_\gamma)$ coincides with $\mathrm{sign}(\kappa^\mathrm{SA})$. 
        \item For any $F\in G^\mathrm{A}$, we have $u_{F\gamma}={\det (DF)} ^{\frac{1}{3}}u_\gamma$, $\sigma_{F\gamma}=\sigma_\gamma$ and $\chi_{F\gamma}=\chi_\gamma$. 
    \end{enumerate}
\end{lemma}
\begin{proof}
    First, 
    %So $\kappa^\mathrm{SA}_\gamma(u)^{-\frac{3}{2}}$ takes value in $\mathbb{R}\cup \imaginary\mathbb{R}$. 
    %As ignoring positive terms in $\chi_\gamma^2 =(\kappa^\mathrm{SA}_\gamma)^{-3}(\kappa^\mathrm{SA}_\gamma)_u^2$ yields just one term $\kappa^\mathrm{SA}_\gamma$, we have (a).
    we have (a) from $\chi_\gamma^2 %=(\kappa^\mathrm{SA}_\gamma)^{-3}(\kappa^\mathrm{SA}_\gamma)_u^2
    =\vlacket{(\kappa^\mathrm{SA}_\gamma)^{-2}{(\kappa^\mathrm{SA}_\gamma)_u}}^2 \kappa^\mathrm{SA}_\gamma$, where $\kappa^\mathrm{SA}_\gamma$ and $(\kappa^\mathrm{SA}_\gamma)_u$ are real-valued.
    Next, for (b), we have
    \begin{align}
        u_{F\gamma}(t)=\int\det ((F\gamma)_t,(F\gamma)_{tt})^\frac{1}{3}dt&=\det (DF)^{\frac{1}{3}}\int\det (\gamma_t,\gamma_{tt})^\frac{1}{3}dt \nonumber
        \\
        &=\det (DF)^{\frac{1}{3}}u_{\gamma}(t).
    \end{align}
    Denote $u_\gamma:=u$, $k:=\det (DF)^\frac{1}{3}$ and $u_{F\gamma}=ku:=v$. 
    By noticing that $u_\gamma$ and $\kappa^\mathrm{SA}_\gamma=-\gamma_{uuu}/\gamma_u$ are invariant under the action of $\det (DF)^{-\frac{1}{2}}F \in G^\mathrm{SA}$, we have
    \begin{align}
        \kappa^\mathrm{SA}_{F\gamma}&=-\frac{(F\gamma)_{vvv}}{(F\gamma)_v}
        =-\frac{k^{-3}\{\det (DF)^{-\frac{1}{2}}F \gamma\}_{uuu}}
        {k^{-1}  \{\det (DF)^{-\frac{1}{2}}F\gamma\}_u} = k^{-2}\kappa^\mathrm{SA}_{\det (DF)^{-\frac{1}{2}}F\gamma}= k^{-2}\kappa^\mathrm{SA}_\gamma,
        \\
        \sigma_{F\gamma}&=\int(\kappa^\mathrm{SA}_{F\gamma})^{\frac{1}{2}}dv=
        \int k^{-1} (\kappa^\mathrm{SA}_{\gamma})^{\frac{1}{2}} (kdu)
        =\sigma_\gamma,
        \\
        \chi_{F\gamma}&=(\kappa^\mathrm{SA}_{F\gamma})^{-\frac{3}{2}}(\kappa^\mathrm{SA}_{F\gamma})_v=(k^{-2}\kappa^\mathrm{SA}_{\gamma})^{-\frac{3}{2}}k^{-1}(k^{-2}\kappa^\mathrm{SA}_{\gamma})_u=\chi_\gamma,
        \end{align}
        which is the conclusion. 
\end{proof}
\subsection{LAC and self-affinities}
In this subsection, we introduce the LAC and its self-affinity as a characterizing property. 
\begin{definition}
     A LAC of slope $\alpha$ is a planar curve defined by
     \begin{align}\label{LACeq}
        \kappa^\mathrm{E}(s)=\left\{\begin{array}{ll}
             (\xi s+\eta)^{-\frac{1}{\alpha}} & (\alpha\neq0), \\[2mm]
             e^{\xi s+\eta}&  (\alpha=0),
        \end{array}\right.
    \end{align}
    for some constants $\xi, \eta\in\mathbb{R}$. 
\end{definition}

\begin{lemma}[\cite{KK}, the Miura self-affinity \cite{Miura2006}]\label{LemmaMSA}
    A curve $\gamma(s)$ is a LAC of slope $\alpha$ if and only if there exists a reparametrization $t=t(s)$  such that for any $t,\varepsilon\in\mathbb{R}$, it follows that
    \begin{align}
    \left\{\begin{array}{ll}
             \kappa^\mathrm{E}(t+\varepsilon)=e^{-\varepsilon} \kappa^\mathrm{E}(t), \\[2mm]
             s_t(t+\varepsilon)=e^{\alpha \varepsilon}s_t(t).
        \end{array}\right.\label{eqMSA}
    \end{align}
\end{lemma}
The above property is called the Miura self-affinity (MSA), which represents an affine relation between two curves $\gamma(t)$ and $\gamma(t+\varepsilon)$ along their Euclidean frames for any shift parameter $\varepsilon\in\mathbb{R}$. 

We note that \eqref{eqMSA} is equivalent to the logarithmic curvature graph \cite{Kanaya2003en, Miura2006} 
\begin{align}
    \Gamma(t)=(X(t),Y(t))= \Big(-\log\kappa^\mathrm{E}(t),\log \left|\frac{ds}{d\log \kappa^\mathrm{E}(t)}\right|\Big),
    \label{LCG}
\end{align} being a line of slope $\alpha$, which is directly verified as follows. 
We have from \eqref{eqMSA} and its $t$-differentiation 
\begin{align}
    \kappa^\mathrm{E}_t(t+\varepsilon)&=e^{-\varepsilon} \kappa^\mathrm{E}_t(t), \\\kappa^\mathrm{E}_s (t+\varepsilon)&=\kappa^\mathrm{E}_t(t+\varepsilon)/s_t(t+\varepsilon)=e^{(\alpha+1)\varepsilon} \kappa^\mathrm{E}_t(t)/s_t(t)=e^{(\alpha+1)\varepsilon}\kappa_s^\mathrm{E}(t), 
    \end{align}
    which yields, by noticing $\log\left| {ds}/{d\log \kappa^\mathrm{E}(t)}\right|=\log\left|{\kappa^\mathrm{E}}/{\kappa^\mathrm{E}_s}\right|$, 
    \begin{align}
       \log \left|\frac{\kappa^\mathrm{E}(t+\varepsilon)}{\kappa_s^\mathrm{E}(t+\varepsilon)}\right|-\log \left | \frac{\kappa^\mathrm{E}(t)}{\kappa_s^\mathrm{E}(t)}\right|=-\varepsilon
        =\alpha\placket{-\log\kappa^\mathrm{E}(t+\varepsilon)+\log\kappa^\mathrm{E}(t)},
    \end{align}
    or 
    \begin{align}
         Y(t+\varepsilon)-Y(t)&=\alpha(X(t+\varepsilon)-X(t)),\label{LCGline}
    \end{align}
    for arbitrary $t,\varepsilon\in\mathbb{R}$. 
    Then, \eqref{LCGline} shows that the LCG is a line of slope $\alpha$. 
\begin{example}\label{exlogspiral}
A \textit{logarithmic spiral} $\gamma(w)=e^{(a+\imaginary b)w}$, $a, b\in{\mathbb{R}}\setminus\{0\}$ is represented in terms of arc length by
    \begin{align}
        \label{eqlogspiral}\gamma(s)=\exp\vlacket{\placket{1+\imaginary \frac{b}{a}}\log\placket{1+\frac{s}{\sqrt{1+(b/a)^2}}}},
    \end{align}
    where $s(w)=\sqrt{1+(b/a)^2}(e^{a
    w}-1)$. 
    Noticing $\gamma_s=e^{\imaginary\theta}$, the turning angle $\theta$ and the curvature $\kappa^\mathrm{E}$ are given by
    \begin{align}
        %\gamma_s(s)&=\exp\vlacket{\imaginary\frac{b}{a}\log\placket{1+\frac{s}{\sqrt{a^2+b^2}}}},\\
        \theta&=\frac{1}{\imaginary}\log \gamma_s=\frac{b}{a}\log\placket{1+\frac{s}{\sqrt{1+(b/a)^2}}}+\mathrm{const.},\label{thetalogspiral}
        \\
        \kappa^\mathrm{E}&=\theta_s=\frac{b/a}{s+\sqrt{1+(b/a)^2}}=:(\xi s+\eta)^{-1}.\label{kappalogspiral}
    \end{align}
    As $\kappa^\mathrm{E}(s)$ is a reciprocal of a linear function, the logarithmic spiral is a LAC of slope $1$.
    Then,    the MSA is demonstrated directly as follows \cite{Miura2006}. We introduce the reparametrization $t=t(s)$ by $s(t)=\frac{\eta}{\xi}(e^{-t}-1)$ so that 
    \begin{equation}
    s_t(t)=-\frac{\eta}{\xi}e^{-t},\quad \kappa^\mathrm{E}(t)=\placket{\xi \frac{\eta}{\xi}(e^{-t}-1)+\eta}^{-1}=\eta e^{t}.
    \end{equation}
    which immediately implies the MSA, 
    $s_t(t+\varepsilon)= e^{-\varepsilon} s_t(t)$ and $\kappa^\mathrm{E}(t+\varepsilon)=e^{\varepsilon}\kappa^\mathrm{E}(t)$.
 %   By introducing the reparametrization $t=t(s)$ so that $s(t)=\frac{\eta}{\xi}(e^{-t}-1)$, we have the MSA as follows \cite{Miura2006}. 
 %   We have $s_t(t)=-\frac{\eta}{\xi}e^{-t}$ and $s_t(t+\varepsilon)= e^{-\varepsilon} s_t(t)$. 
%    Furthermore, \eqref{kappalogspiral} is represented by 
%    \begin{align}
%        \kappa^\mathrm{E}(t)=\placket{\xi \frac{\eta}{\xi}(e^{-t}-1)+\eta}^{-1}=\eta e^{t},
%    \end{align}
%    for which we have $\kappa^\mathrm{E}(t+\varepsilon)=e^{\varepsilon}\kappa^\mathrm{E}(t)$. 
\end{example}

In the early stage of the study of the LAC, Harada et.\,al.\,\cite{Harada1999, Harada1995} claimed that the LAC possesses another kind of self-affinity %that any subcurve is affinely equivalent to the whole curve. 
that any subcurve, with arbitrary scale change in the horizontal and the vertical directions, is affinely equivalent to the whole curve. 
Noticing the former scale change belongs to $G^\mathrm{A}$, their self-affinity is defined rigorously as follows. 
\begin{definition}[the Harada self-affinity \cite{Harada1999, Harada1995}]
        A curve $\gamma(s):\mathbb{R}\rightarrow {\mathbb{C}}$ possess the \textit{Harada self-affinity (HSA)} 
     if its arbitrary open subcurve is affinely equivalent to the whole curve. 
     That is, 
     for any open subinterval $I\subset \mathbb{R}$, there exists a reparametrization $t=t(s):\mathbb{R}\rightarrow I$ and an affine map $F_I:{\mathbb{C}}\rightarrow {\mathbb{C}}$ such that $\gamma(s)=F_I\,\gamma(t(s))$ on $\mathbb{R}$. 
\end{definition}
Though the HSA was found not to be the self-affinity of the LAC, the authors proved the following proposition. 
\begin{proposition}[\cite{KK}]\label{HSAparabola}
    A curve possesses the HSA if and only if it is either a line %, a circle, 
     or a parabola. 
\end{proposition}
Then, in view of Lemma \ref{quad}, Proposition \ref{HSAparabola} suggests that the HSA should be discussed in the framework of equiaffine geometry. 
Accordingly, the HSA can be also relaxed into a \textit{self-affinity} of constant curvature curves, namely, quadratic curves.
To this end, let us demonstrate that the logarithmic spiral possesses another self-affinity. 
The equation \eqref{eqlogspiral} is rewritten in terms of the turning angle \eqref{thetalogspiral} by
    \begin{align}
        %s(\theta)&=\sqrt{a^2+b^2}(e^{(a/b)\theta}-1), \\
        \gamma(\theta)&=e^{{(\frac{a}{b}+\imaginary) }(\theta+\mathrm{const.})}.  
        %\exp\vlacket{\placket{1+\imaginary \frac{b}{a}}t}. 
    \end{align}
    Then, for any $\theta,\varepsilon\in\mathbb{R}$, we observe
    \begin{align}
        \gamma(\theta+\varepsilon)
        &=e^{(\frac{a}{b}+\imaginary)\varepsilon}\gamma(\theta)
        %&=\exp\vlacket{\placket{\frac{a}{b}+\imaginary}\varepsilon}\gamma(\theta)
        =e^{\frac{a}{b}\varepsilon}\begin{pmatrix} \cos\varepsilon&-\sin\varepsilon\\\sin\varepsilon&\cos\varepsilon\end{pmatrix}\gamma(\theta), \label{ESAlogspiral}
    \end{align}
%It is formulated by the property that the original curve is affine equivalent to its arbitrary shift, just like the logarithmic spiral \eqref{ESAlogspiral}, 
which is a self-affinity in the sense that the two curves $\gamma(\theta +\varepsilon)$ and $\gamma(\theta)$ are related by a similarity transformation depending on the shift parameter $\varepsilon$. 
In the following, we define this self-affinity in a more general setting.  
\if0
\begin{definition}
    A \textit{Lie group} is a group $G$ equipped with a structure of differentiable manifold on $G$ such that the mapping $G\times G:(x,y)\mapsto x^{-1}y$ is smooth. 
    A \textit{Lie algebra} of $G$ is a tangent space $T_gG$, $g\in G$. 
\end{definition}
\fi
\begin{definition}\label{defESA}
    A curve $\gamma(w):\mathbb{R}\rightarrow {\mathbb{C}}$ possesses the \textit{extendable self-affinity (the ESA)} with respect to a subgroup $G$ of the Lie group $G^\mathrm{A}$ \eqref{defGA} if there exists a reparametrization $w=w(t)$ 
    % $t(w)$ of $\gamma(w)$ 
    and  
    a differentiable map $F(\varepsilon):\mathbb{R}\rightarrow G$ such that $F(0)=\mathrm{id}$ and 
\begin{equation}\label{eqESA}
    \gamma(t+\varepsilon)=F(\varepsilon)\,\gamma(t),
\end{equation}\label{ESA}
 for any $t,\varepsilon\in\mathbb{R}$. 
The parameter $t$ is referred to as the \textit{ESA-parameter of $\gamma$}.
\end{definition}
We note that the group $G$ can be chosen as a more general Lie group, namely 
the group acting on curves arising from affine action on their Frenet frames  \cite[(4.2)]{KK}, 
%the group of the pointwise correspondence \cite[(4.2)]{KK} arising from the MSA \eqref{eqMSA} which is correspondence of charasteristics, 
%between $\gamma(t+\varepsilon)$ and $\gamma(t)$ for a LAC $\gamma(t)$, 
but we do not consider such a case in this paper. 
\begin{proposition}[\cite{KK}]\label{ESAquad}
    A curve $\gamma(t):\mathbb{R}\rightarrow {\mathbb{C}}$ possesses the ESA with respect to the group $G^\mathrm{SA}$ and with the ESA-parameter being the equiaffine parameter $u(t)$ if and only if its equiaffine curvature is constant, equivalently it is a quadratic curve. 
\end{proposition}

\section{Curves with ESA in terms of equiaffine geometry}
\subsection{Concrete form of the ESA-parameter and the coefficient matrix}
In the following, we extend Proposition \ref{ESAquad} without any restriction on the ESA-parameter $t$ and show that a specific relation between $t$ and $u$ is identified.
%revise
For independent real variables $t,\varepsilon,\varepsilon'$, using the ESA \eqref{eqESA}, $\gamma(t+\varepsilon+\varepsilon')$ and its derivatives are represented by
\begin{align}
    F(\varepsilon)F(\varepsilon')\gamma(t)&=F(\varepsilon')F(\varepsilon)\gamma(t)=F(\varepsilon + \varepsilon')\gamma(t),
    \\
    DF(\varepsilon)DF(\varepsilon')\gamma_t(t)&=DF(\varepsilon')DF(\varepsilon)\gamma_t(t)=DF(\varepsilon + \varepsilon')\gamma_t(t),\label{DFSA1}
    \\DF(\varepsilon)DF(\varepsilon')\gamma_{tt}(t)&=DF(\varepsilon')DF(\varepsilon)\gamma_{tt}(t)=DF(\varepsilon + \varepsilon')\gamma_{tt}(t). \label{DFSA2}
\end{align}
We recall that $\gamma(t)$ is non-degenerate $\det(\gamma_t,\gamma_{tt})\neq0$ and $DF(\varepsilon)$ takes value in $ \mathrm{M}(2,\mathbb{R})$. 
Then, \eqref{DFSA1} and \eqref{DFSA2} yield the additive law $DF(\varepsilon+\varepsilon')=DF(\varepsilon)DF(\varepsilon')=DF(\varepsilon')DF(\varepsilon)$, whose $\varepsilon'$-differentiation at $\varepsilon'=0$ is  
\begin{align}
    \frac{d}{d\varepsilon'}DF(\varepsilon+\varepsilon')\Big|_{\varepsilon'=0}=DF_\varepsilon(\varepsilon)=DF_\varepsilon(0)DF(\varepsilon),
\end{align}
solving which we identify $DF(\varepsilon)=e^{\varepsilon A}$ where $A:=DF_\varepsilon(0)\in \mathrm{M}(2,\mathbb{R})$. %\mathrm{M}(2,\mathbb{R})$. 
Then, 
from \eqref{DFSA1} and \eqref{DFSA2} we have 
\begin{align}\label{SAframe}
    \Phi(t+\varepsilon)=e^{\varepsilon A}\Phi(t), \quad \Phi(t):=(\gamma_{t}(t), \gamma_{tt}(t)). 
    %\\
    %\Phi:=(\gamma_{w},\gamma_{ww})= (\gamma_{u},\gamma_{uu})\begin{pmatrix}        w_u& w_{uu}\\        0& w_{u}    \end{pmatrix}=:\Phi^\mathrm{SA}(w)U(w),
\end{align}
Taking $\varepsilon$-differentiation of \eqref{SAframe} at $\varepsilon=0$ yields
$\Phi_t=A\Phi$, from which we have $\Phi(t)=e^{tA}\Phi(0)$. 
Then, solving $u_t(t)=\det\Phi(t)^\frac{1}{3}=\det\Phi(0)^\frac{1}{3}e^{\frac{t}{3} (\mathrm{trA})}$ yields
\begin{align}
    u(t)=\int ke^{(\mathrm{trA})t}dt=\begin{cases}
        kt + \mathrm{const.} & \text{if $\mathrm{tr}A=0$}
        \\
        \frac{k}{\mathrm{tr}A}e^{t(\mathrm{trA})} + \mathrm{const.} & \text{if $\mathrm{tr}A\neq 0$},
    \end{cases} \quad k=\det\Phi(0)^{\frac{1}{3}}.\label{uESA}
\end{align}
%We note that $\mathrm{tr}A=0$ if and only if $DF(\varepsilon)$ belongs to the Lie algebra $\mathfrak{sl}(2,\mathbb{R})$ of the Lie group $\mathrm{SL}(2,\mathbb{R})$. 
Noticing that $\mathrm{tr} A = 0$ holds if and only if $DF(\varepsilon)$ belongs to the Lie algebra $\mathfrak{sl}(2,\mathbb{R})$ of the Lie group $\mathrm{SL}(2,\mathbb{R})$, we obtain the following result.
%We have obtained the following result. 
\begin{lemma}\label{ESASA}
    Let $\gamma(t)$ be a non-degenerate curve with the ESA $\gamma(t+\varepsilon)=F(\varepsilon)\gamma(t)$, and $F(\varepsilon)$ act on $\gamma_t$ by affine action represented by $DF(\varepsilon)\in \mathrm{GL}(2,\mathbb{R})$. 
    Then, up to affine transformations $u\mapsto ku+l$, $t\mapsto k't+l'$, the equiaffine parameter $u(t)$ coincides with either the ESA-parameter $t$ or its exponential function $e^t$. The former occurs if and only if $F(\varepsilon)$ is $G^\mathrm{SA}$-valued.  
\end{lemma}
Lemma \ref{ESASA} implies that it is natural to %extend $G=G^\mathrm{A}$. 
consider the ESA in $G^\mathrm{A}$, namely, in affine geometry, which will be discussed in the next subsection.
\subsection{ESA in affine geometry and statement of main theorem}

By Lemma \ref{affinv} (b), for a curve with the ESA in $G^\mathrm{A}$ and $t,\varepsilon\in \mathbb{R}$, we have 
\begin{align}
    \chi_\gamma(t+\varepsilon)=\chi_{F(\varepsilon)\gamma}(t)=\chi_\gamma(t),\label{SAchi}
\end{align}
which implies the constantness of the affine curvature $\chi_\gamma(t)$.
More strongly, we claim that the ESA in $G^\mathrm{A}$ uniquely characterizes the curves with constant affine curvature as follows. 
\begin{theorem}\label{main}
    A curve possesses the ESA with respect to the affine transformation group $G^\mathrm{A}$ if and only if it is up to $G^\mathrm{A}$ either of
    \begin{enumerate}
        \item the graph of a power function $y=x^p$, 
        \item the graph of $y=\log x$, 
        \item the graph of $y=x\log x$,
        \item a logarithmic spiral, or
        \item a quadratic curve.
    \end{enumerate}
\end{theorem}
We note that the graph of $y=e^x$ is affinely equivalent to the graph of $y=\log x$, which belongs to our class. 
The family of curves in Theorem \ref{main} may be regarded as an alternate family of \textit{``aesthetic curves"}, since it has self-affinity and contains logarithmic spiral which is a special case of the LAC. 
Furthermore, this family contains the quadratic curve which play a fundamental role in CAGD. 
Figure \ref{fig:fig1} illustrates some examples of the curves in the class. 
\if0
\begin{figure}[h]
\centering
    \begin{minipage}[b]{0.49\columnwidth}
        \centering
       \includegraphics[height=0.6\columnwidth]{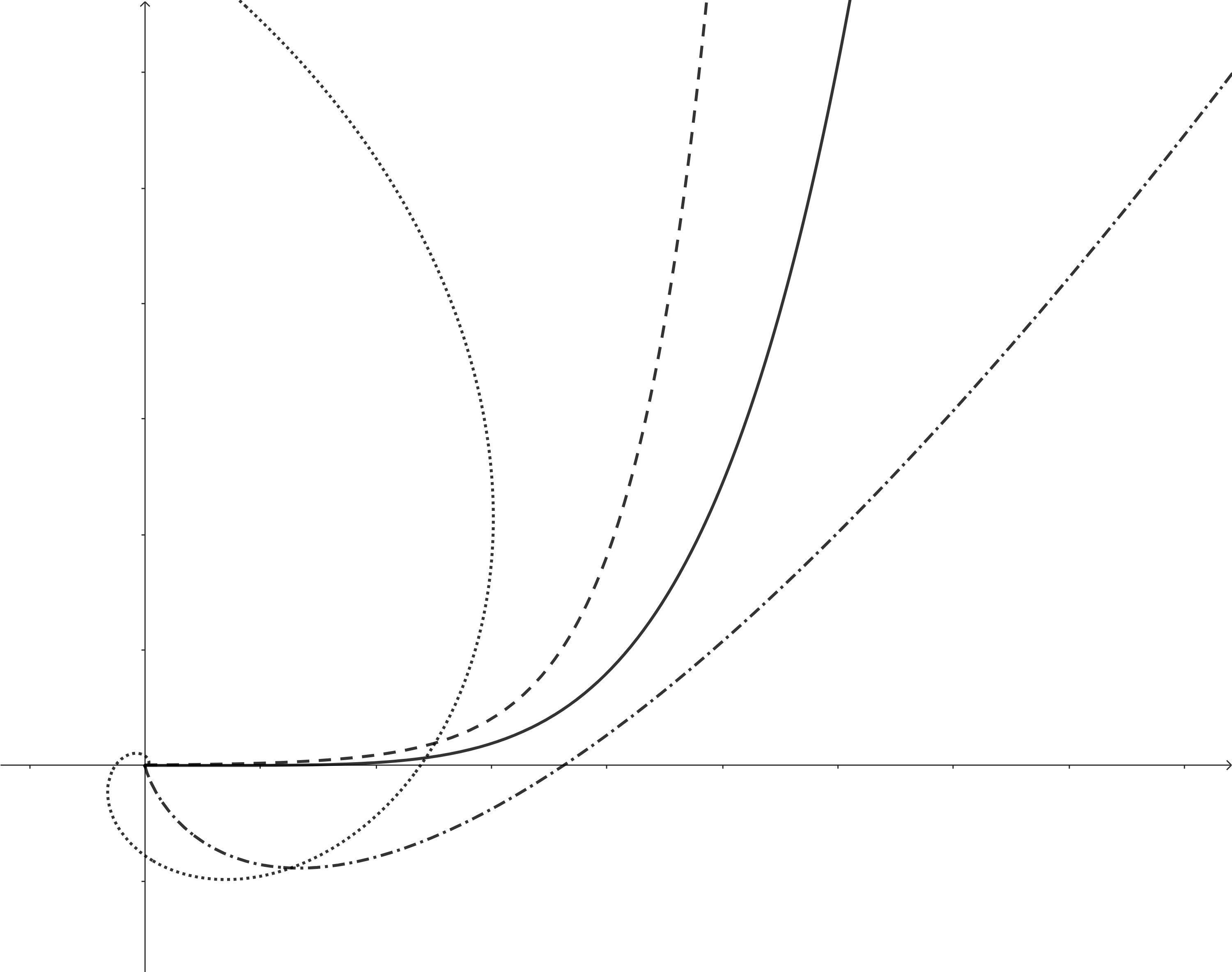}  
    \end{minipage}
    \begin{minipage}[b]{0.49\columnwidth}
        \centering
       \includegraphics[height=0.6\columnwidth]{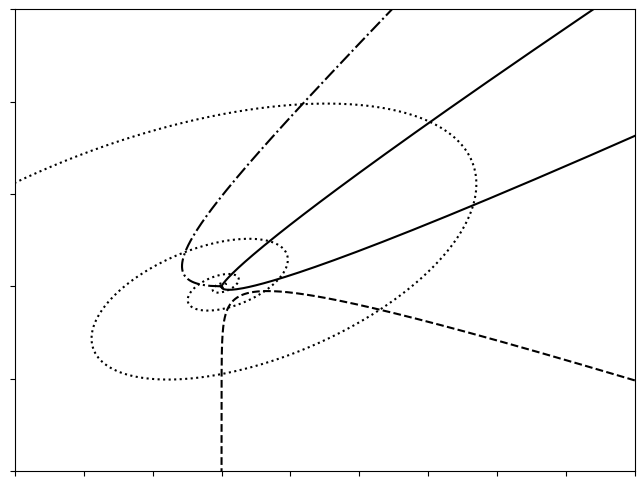}
    \end{minipage}
       \caption{family of curves determined by the ESA in $G^\mathrm{A}$. Solid line: power function, dotted line: logarithmic spiral, dashed line: exponential (logarithmic) function, dashdotted line: $x\log x$. 
       Left: Comparison of curves starting from the origin. Right: Affine deformation applied to each curve.}
   \label{fig:fig1}
\end{figure}
\fi
\begin{figure}[h]
\centering
       \includegraphics[width=\columnwidth]{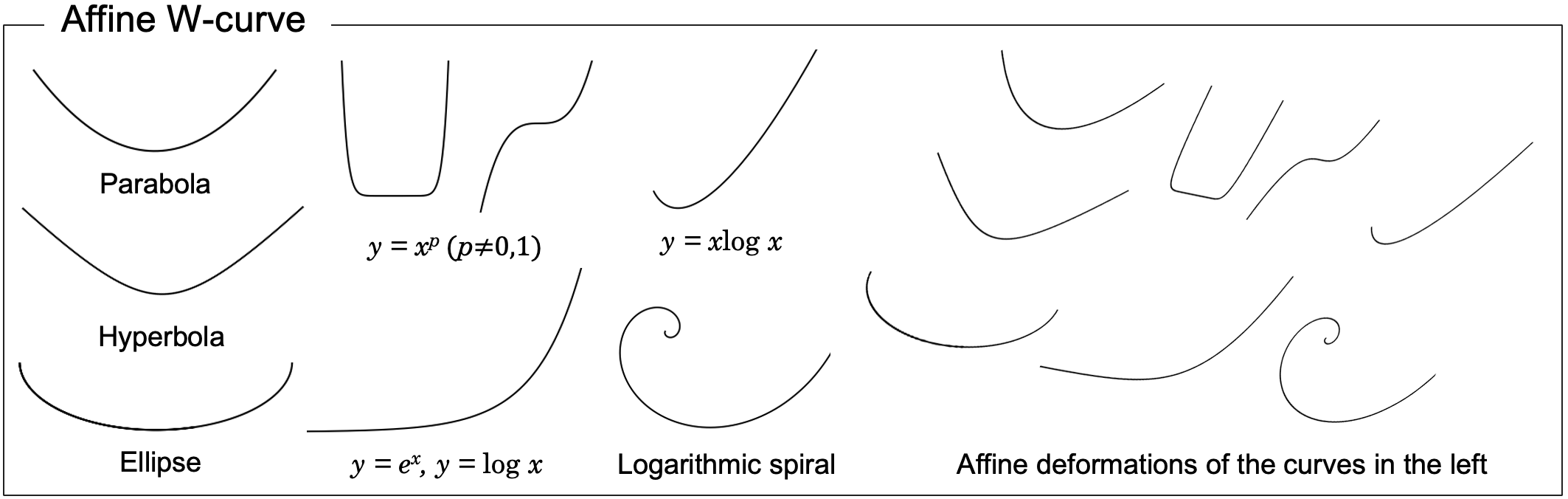}
       \caption{affine W-curves: family of curves determined by the ESA in $G^\mathrm{A}$.}
   \label{fig:fig1}
\end{figure}
The family shown in Theorem \ref{main} have been pointed out by Blaschke \cite{Blaschke} as a special case of so-called the Klein-Lie \textit{W-curve} \cite{LieKlein} in projective differential geometry. 
We call the curves in Theorem \ref{main} the \textit{affine version of Klein-Lie W-curves} or the \textit{affine W-curves} in short in this paper.
We carefully discuss an identification of the family \cite{Blaschke}, which is necessary for a direct proof of the ESA in $G^\mathrm{A}$. 

We prove Theorem \ref{main} by showing the following proposition. 

\begin{proposition}\label{keyProp}
    For a non-degenerate planar curve $\gamma(t)$, the following are equivalent. 
    \begin{enumerate}
        \item The affine curvature $\chi(t)$ of $\gamma(t)$ is constant.  
        \item The equiaffine curvature of $\gamma(t)$ is of the form $\kappa^\mathrm{SA}(u)=\delta(\xi u+\eta )^{-2}$,  $\delta\in\{\pm 1\}$, $\xi,\eta\in\mathbb{R}$.
        \item $\gamma(t)$ belongs to the affine W-curve. 
        \item $\gamma(t)$ posesses the ESA in $G^\mathrm{A}$. 
    \end{enumerate}
    As a trivial case, $\chi=0$ in (a), $\kappa^\mathrm{SA}$ being constant in (b), $\gamma$ being quadratic curve in (c), and the ESA formula of the quadratic curve \cite[Theorem 4.2]{KK} are equivalent. 
\end{proposition}

\subsection{Curves with constant affine curvature}
In this and next subsections, we prove Proposition \ref{keyProp}. 
We first note that (d)$\Rightarrow$(a) follows from the discussion on \eqref{SAchi}. 
We then discuss the parts (a)$\Rightarrow $(b)$\Rightarrow$(c) in Proposition \ref{keyProp}.

    Let $\chi=(\kappa^\mathrm{SA})^{-\frac{3}{2}}\kappa^\mathrm{SA}_u$ be arbitrarily constant in $(\mathbb{R}\cup \imaginary\mathbb{R})\setminus\{0\}$.
    By Lemma \ref{affinv} (a), $\chi\neq0$ and smoothness of $\gamma$, we have  
    $\delta:=\mathrm{sign}(\kappa^\mathrm{SA}(u))\in\{\pm1\}$  %,     $\kappa^\mathrm{SA}_u\neq0$ anywhere, 
    constant with ${\chi}/{\sqrt{\delta}}\in\mathbb{R}$. 
    Then, solving the differential equation $(\kappa^\mathrm{SA})^{-\frac{3}{2}}\kappa^\mathrm{SA}_u=\chi$ (constant) in a standard manner yields, 
    \begin{align}
        %\int^u_{u_0} (\kappa^\mathrm{SA})^{-\frac{3}{2}}d\kappa^\mathrm{SA}&=\chi\int^u_{u_0} du, \\
        %-2(\kappa^\mathrm{SA})^{-\frac{1}{2}}+2(\kappa^\mathrm{SA}_0)^{-\frac{1}{2}}&=\chi (u-u_0), \\
        \kappa^\mathrm{SA}={\vlacket{\chi(u-u_0)-2(\kappa^\mathrm{SA}_0)^{-\frac{1}{2}}}^{-2}}&={\delta\vlacket{\frac{\chi}{\sqrt{\delta}}(u-u_0)-2|\kappa^\mathrm{SA}_0|^{-\frac{1}{2}}}}^{-2},\label{kappaSA-2}
    \end{align}
    where $u_0, \kappa^\mathrm{SA}_0$ are arbitrary constants. 
    Here we have used $\mathrm{sign}(\chi^2)=\mathrm{sign}(\kappa^\mathrm{SA})=\delta$, ${\chi}/{\sqrt{\delta}}\in\mathbb{R}$ and $\delta^2=1$.
    Then, \eqref{kappaSA-2} is represented by the form $\delta(\xi u+\eta)^{-2}$. 
    Thus, we proved (a)$\Rightarrow$(b). 

We next prove (b)$\Rightarrow$(c) by showing that $\kappa^\mathrm{SA}(u)=\delta(\xi u+\eta )^{-2}$ yields the affine W-curve. 
We may assume $\xi\neq0$ to avoid the trivial case in view of Lemma \ref{quad}. 
By using the translation ambiguity of the equiaffine parameter $u$, we can choose $\eta=0$ without loss of generality. 
By Lemma \ref{lemma:repformula}, we have Euler's differential equation
\begin{align}\label{Euler}
    z_{uu}+\delta(\xi u)^{-2}z=0.
\end{align}
Solving \eqref{Euler} in a standard manner yields 
\begin{align}
    &z(u)=\exp\vlacket{\log u\,\frac{\xi\pm\sqrt{\xi^2-4\delta}}{2\xi}  }=:u^{\frac{1}{2}\pm\frac{\omega}{2}},\quad \omega=\sqrt{\Delta}:=\sqrt{1- \frac{4\delta}{\xi^2}}.
\end{align}
We note that $\omega$ takes arbitrary value in $\mathbb{R}\setminus \{1\}$ as $\delta\in\{\pm1\}$ and $\xi\in\mathbb{R}\setminus\{0\}$ vary. 
According to the sign of the discriminant $\Delta$, a classification of basis of real-valued solution of \eqref{Euler} is given by
\begin{align}
    &(f(u),g(u))=\begin{cases}
        (u^{\frac{1}{2}+\frac{\omega}{2}},u^{\frac{1}{2}-\frac{\omega}{2}})&(\Delta>0),
        \\[2mm]
        \placket{\sqrt{u},\sqrt{u}\log u}&(\Delta=0),
        \\[2mm]
        %\placket{\sqrt{u}\cos\placket{\frac{\omega}{2\imaginary}\log u},\sqrt{u}\sin\placket{\frac{\omega}{2\imaginary}\log u}}
        \placket{\mathrm{Re}\placket{u^{\frac{1}{2}+\imaginary\frac{|\omega|}{2}}},\mathrm{Im}\placket{u^{\frac{1}{2}+\imaginary\frac{|\omega|}{2}}}}&(\Delta<0),\label{sLACcurvescoord}
    \end{cases}
\end{align}
where $\det\begin{pmatrix}
    f&f_u\\g& g_u
\end{pmatrix}\neq 0$. 
%Up to translations, we may assume $\gamma(u)=x(u)+\imaginary y(u)$  $x(u_0)=y(u_0)=0$ without loss of generality. 
Then, from Lemma \ref{lemma:repformula}, integrating \eqref{sLACcurvescoord} by $u$ yields
\begin{align}
    (x(u),y(u))&= \begin{cases}
        \blacket{\dfrac{u^{\frac{3}{2}+\frac{\omega}{2}}}{\frac{3}{2}+\frac{\omega}{2}},\, \dfrac{u^{\frac{3}{2}-\frac{\omega}{2}}}{\frac{3}{2}-\frac{\omega}{2}}}^u_{u_0}&(\Delta>0, \omega\neq \pm3),
        \\[4mm]
        \blacket{\frac{1}{3}u^{3},\log u}^u_{u_0}&(\Delta>0, \omega=\pm3),
        \\[2mm]
         \blacket{\frac{2}{3}u^\frac{3}{2}, \frac{2}{3}u^\frac{3}{2}\placket{\log u-\frac{2}{3}}}^u_{u_0}&(\Delta=0),
        \\[2mm]
        \blacket{\mathrm{Re}\placket{\dfrac{u^{\frac{3}{2}+\imaginary\frac{|\omega|}{2}}}{{\frac{3}{2}+\imaginary\frac{|\omega|}{2}}}},\mathrm{Im}\placket{\dfrac{u^{\frac{3}{2}+\imaginary\frac{|\omega|}{2}}}{{\frac{3}{2}+\imaginary\frac{|\omega|}{2}}}}}^u_{u_0}&(\Delta<0),
    \end{cases}\label{sLACcurves}
\end{align}
 where $\gamma(u)=x(u)+\imaginary y(u)$ and $\gamma(u_0)=0$, up to affine transformations. 
 Thus, $\gamma$ is classified into the graph of either of 
\begin{align}
    \begin{cases}
        y=x^p 
        &(\Delta>0, \omega\neq \pm3),\\
        y=\log x&(\Delta>0, \omega=\pm3),\\
        y=x\log x&(\Delta=0),
        \\
        x(t)+\imaginary y(t)=e^{(a+\imaginary b)t
        }&(\Delta<0),\; a,b\in\mathbb{R},
    \end{cases}\label{sLACcurves2}
\end{align}
%Figures 1-4 show the example of the above curves. 
where any exponent $p\in\mathbb{R}\setminus\{0,-1,2\}$ is realized by $\omega=3\frac{1-p}{1+p}\in\mathbb{R}\setminus \{1,\pm3\}$. 
%\\
%\note{$\alpha\rightarrow (\infty, 0, -1,1,2)$ corresponds to $\omega\rightarrow (-3,3,\infty,0, 1)$}
We note that the cases $p=-1$ (the hyperbola), $p=1$ (the line), $p=2$ (the parabola) are degenerate but captured by the quadratic curve and the line. %(Proposition \ref{quad} and Proposition \ref{HSAparabola}). 
So any possibilities in \eqref{sLACcurves2} appears as a solution to the constantness of affine curvature $\chi$.
Now we have proved (b)$\Rightarrow$(c).

We also note that the fourth case in \eqref{sLACcurves2} is a logarithmic spiral as shown in Example \ref{exlogspiral}, and it is a circle (ellipse) when $a=0$. 

\subsection{Explicit formula of the ESA in $G^\mathrm{A}$}

\if0
Solving \eqref{Euler} in a standard manner yields 
\begin{align}
    &z(u)=\exp\vlacket{\log u\,\frac{\xi\pm\sqrt{\xi^2+4\varepsilon}}{2\xi}  }=:u^{\frac{1}{2}\pm\omega},
\end{align}
where $\omega=\sqrt{\Big|1+\dfrac{\varepsilon}{4\xi^{2}}\Big|}\geq0$. 
According to sign of the discriminant $\Delta=\xi^2+4\varepsilon$, basis of real-valued solution of \eqref{Euler} is classified by 
\begin{align}
    &(f(u),g(u))=\begin{cases}
        (u^{\frac{1}{2}+\omega},u^{\frac{1}{2}-\omega})&(\Delta>0),
        \\
        \placket{\sqrt{u},\sqrt{u}\log u}&(\Delta=0),
        \\
        \placket{\sqrt{\frac{u}{\omega}}\cos\placket{\omega\log u},\sqrt{\frac{u}{\omega}}\sin\placket{\omega\log u}}&(\Delta<0),\label{sLACcurvescoord}
    \end{cases}
\end{align}
where $\det\begin{pmatrix}
    f&f_u\\g& g_u
\end{pmatrix}\neq 0$. 
We note that as long as $\varepsilon=-1$, any case of $\Delta>0$ with arbitrary exponent $\alpha'=\log_g(f)\in\mathbb{R}^\times$,  $\Delta=0$, and $\Delta<0$ is possible up to affine transformation. 
On the other hand, $\varepsilon=+1$ implies that we have $\Delta>0$ and $\omega>1$ for arbitrary $\xi$. 
\medskip
\fi

In this section, we establish (c)$\Rightarrow$(d) in Proposition \ref{keyProp}, by 
directly verifying that the curve described in \eqref{sLACcurves} possesses the ESA in $G^\mathrm{A}$. 
We continue settings in previous section.
First, up to $G^\mathrm{A}$ the representation \eqref{sLACcurvescoord} is rewritten by
\begin{align}
        (x(u),y(u))&= \begin{cases}
        \ \placket{{u^{\frac{3}{2}+\frac{\omega}{2}}},\, {u^{\frac{3}{2}-\frac{\omega}{2}}}}&(\Delta>0, \omega\neq \pm3),
        \\[2mm]
        \ \placket{u^{3},\log u}&(\Delta>0, \omega=\pm3),
        \\[2mm]
         \ \placket{u^\frac{3}{2}, u^\frac{3}{2}{\log u}}&(\Delta=0),
        \\[2mm]
        \ \placket{\mathrm{Re}\placket{{u^{\frac{3}{2}+\imaginary\frac{|\omega|}{2}}}},\mathrm{Im}\placket{{u^{\frac{3}{2}+\imaginary\frac{|\omega|}{2}}}}}&(\Delta<0).
    \end{cases}\label{sLACcurvescoordsimple}
\end{align}
By Lemma \ref{ESASA}, without the trivial case, the ESA-parameter $t$ is up to translation given by $u(t)=e^{kt+l}$, $k,l\in\mathbb{R}$.
Substituting $u(t)=e^{kt+l}$ to \eqref{sLACcurvescoordsimple} yields
%\note{in progress}
\begin{align}
    \begin{pmatrix}
            x(t)\\[2mm]y(t)
        \end{pmatrix}&=\begin{cases}
        \ \begin{pmatrix}
            e^{(\frac{3}{2}+\frac{\omega}{2})(kt+l)}\\e^{(\frac{3}{2}-\frac{\omega}{2})(kt+l)}
        \end{pmatrix} &(|\Delta|>0,\omega\neq3),%\text{ if }|\xi|=2,
        \\
        \ \begin{pmatrix}
            e^{3(kt+l)}\\kt+l
        \end{pmatrix} &(|\Delta|>0,\omega=3),%\text{ if }|\xi|=2,
        \\
        \ \begin{pmatrix}
            e^{\frac{3}{2}(kt+l)}
            \\(kt+l)e^{\frac{3}{2}(kt+l)}
        \end{pmatrix} &(|\xi|=2),
        \\
        \ \begin{pmatrix}
            e^{\frac{3}{2}(kt+l)}\cos{\frac{|\omega|}{2} (kt+l)}\\
            e^{\frac{3}{2}(kt+l)}\sin{\frac{|\omega|}{2} (kt+l)}
        \end{pmatrix}&(|\xi|<2),%\label{sLACcurvest}
    \end{cases}
\end{align}
from which we have 
\begin{align}
    \begin{pmatrix}
            x(t+\varepsilon)\\[2mm]y(t+\varepsilon)
        \end{pmatrix}&
    =\begin{cases}
        \ e^{\frac{3k\varepsilon}{2}}\begin{pmatrix}
            e^{\frac{k\omega}{2}\varepsilon}&0\\0&e^{-\frac{k\omega}{2}\varepsilon}
        \end{pmatrix}\begin{pmatrix}
            e^{(\frac{3}{2}+\frac{\omega}{2})(kt+l)}\\e^{(\frac{3}{2}-\frac{\omega}{2})(kt+l)}
        \end{pmatrix}&(|\Delta|>0,\omega\neq3),%\text{ if }|\xi|=2,
        \\\ e^{\frac{3k\varepsilon}{2}}
        \begin{pmatrix}
            e^{\frac{3k\varepsilon}{2}}&0\\0&e^{-\frac{3k\varepsilon}{2}}
        \end{pmatrix}\begin{pmatrix}
            e^{{3}(kt+l)}
            \\kt+l
        \end{pmatrix}
        +\begin{pmatrix}
            0\\[2mm]k\varepsilon
        \end{pmatrix}
        &(|\Delta|>0,\omega=3),%\text{ if }|\xi|=2,
        \\
        \ e^{\frac{3k\varepsilon}{2}}
        \begin{pmatrix}
            1&0\\k\varepsilon&1
        \end{pmatrix}\begin{pmatrix}
            e^{\frac{3}{2}(kt+l)}
            \\(kt+l)e^{\frac{3}{2}(kt+l)}
        \end{pmatrix} &(|\xi|=2),
        \\
        \ e^{\frac{3k\varepsilon}{2}}
        \begin{pmatrix}
            \cos \frac{k|\omega| \varepsilon}{2}&-\sin \frac{k|\omega| \varepsilon}{2}\\\sin \frac{k|\omega| \varepsilon}{2}
            &\cos\frac{k|\omega| \varepsilon}{2}
        \end{pmatrix}
        \begin{pmatrix}
            e^{\frac{3}{2}(kt+l)}\cos{\frac{|\omega|}{2} (kt+l)}\\
            e^{\frac{3}{2}(kt+l)}\sin{\frac{|\omega|}{2} (kt+l)}
        \end{pmatrix}&(|\xi|<2).%\label{sLACcurvest}
    \end{cases}\label{sLACcurvest}
\end{align}
The right hand side of \eqref{sLACcurvest} is always represented by the form $F(\varepsilon)\begin{pmatrix}
            x(t)\\[2mm]y(t)
        \end{pmatrix}$, where $F(\varepsilon)\in G^\mathrm{A}$. 
Then, we have shown that the curve described in \eqref{sLACcurves} possesses the ESA in $G^\mathrm{A}$ explicitly.  
This completes the proof of Proposition \ref{keyProp} and Theorem \ref{main}.

We may remark that from $\kappa^\mathrm{SA}(u)=\delta(\xi u+\eta)^{-2}$ the affine arc length \eqref{sigma} is calculated as
\begin{align}
    \sigma(u)=\int \vlacket{\delta(\xi u+\eta)^{-2}}^\frac{1}{2}du =\frac{\sqrt{\delta}}{\xi}\log (\xi u+\eta)+(\mathrm{const.}).
\end{align}
Thus, the ESA-parameter $t$ is identified with the affine arc length $\sigma$ up to multiplication by $\imaginary$ for the affine W-curves, except quadratic curves ($\xi=0$) for which $t$ is identified with the equiaffine arc length $u$. 

The above discussions may imply that the ESA in a Lie group $G$ gives rise to curves  of constant curvature in the Klein geometry associated with $G$, where the ESA-parameter is an invariant parameter under $G$. 
As long as considering self-affinity of planar curves, the ESA in affine geometry and the family of curves that we have identified the affine W-curve, shall be the most general ones. 

\if0
\note{xlogx}
\begin{figure}[h]
\centering
\begin{minipage}[b]{0.49\columnwidth}
    \centering
    \includegraphics[width=0.9\columnwidth]{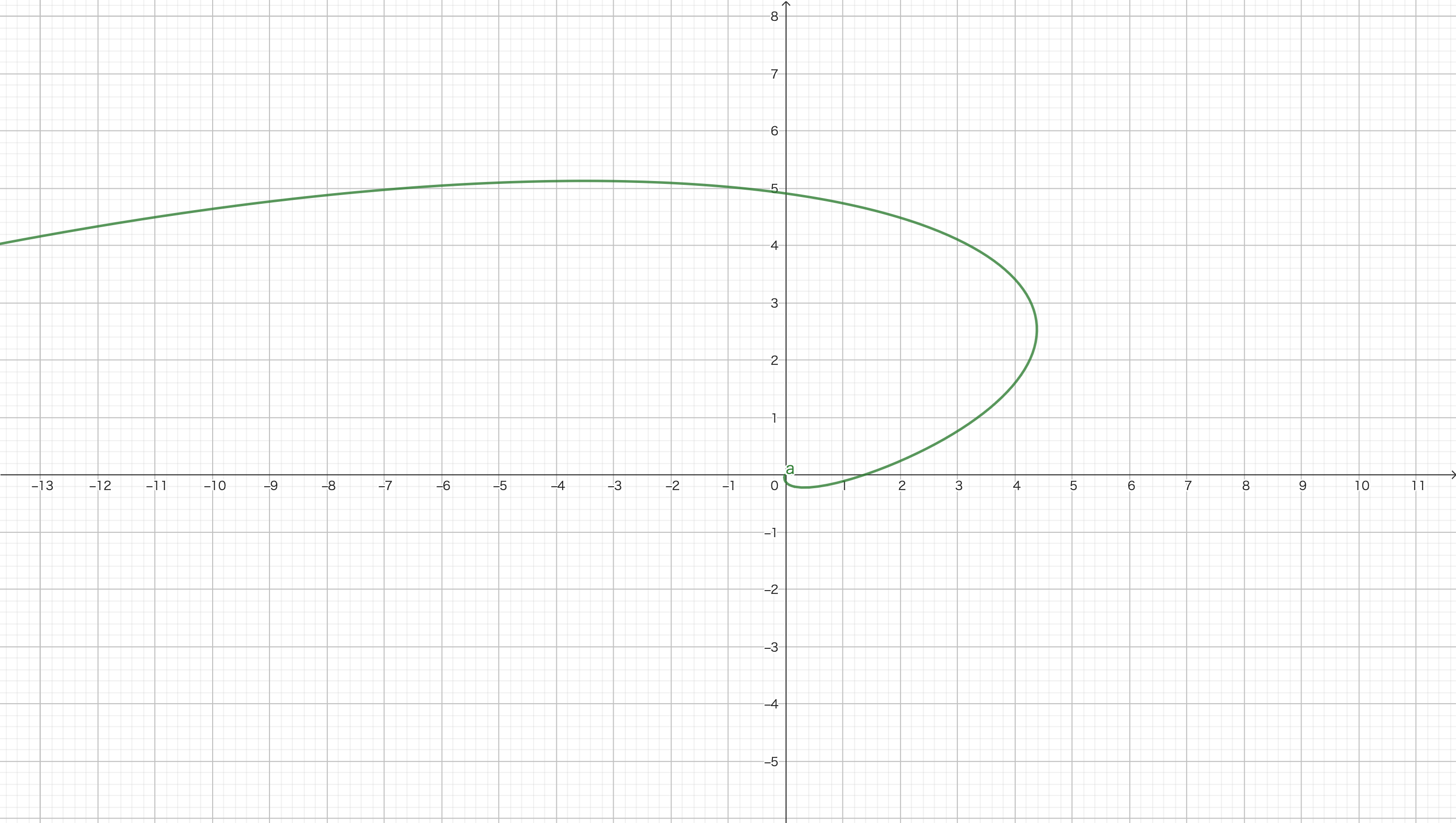}
    \caption{Curve with $\xi=1$.}
    \label{fig:a}
\end{minipage}
\begin{minipage}[b]{0.49\columnwidth}
    \centering
    \includegraphics[width=0.9\columnwidth]{SALAC_real_a_-4_1(1).png}
    \caption{Curve with $\frac{1}{2}\pm\omega=-4,1$. }
    \label{fig:b}
\end{minipage}
\end{figure}
\begin{figure}[h]
\centering
\begin{minipage}[b]{0.49\columnwidth}
    \centering
    \includegraphics[width=0.9\columnwidth]{SALAC_real_b_-2_1(1).png}
    \caption{Curve with $\frac{1}{2}\pm\omega=-2,1$.}
    \label{fig:c}
\end{minipage}
\begin{minipage}[b]{0.49\columnwidth}
    \centering
    \includegraphics[width=0.9\columnwidth]{SALAC_real_c_2_1(1).png}
    \caption{Curve with $\frac{1}{2}\pm\omega=1,2$. }
    \label{fig:d}
\end{minipage}
\end{figure}
\fi

\section{Logarithmic curvature graphs of curves in the class}
In this section, we discuss the relationship between the affine W-curves and Harada-Miura's ``aesthetic curves" in industrial design. 
%Harada et al.\cite{Harada1999} observed that the professional car designers and CAD experts regard as ``aesthetic" have a common property; almost linearity of LCH. 
%They call such curves the {monotonous rhythm curves}. 
%Then, the \textit{exact linearity of the LCG} leads to the definition of LAC. 
%
As discussed following Lemma \ref{LemmaMSA}, the LAC is defined by the \textit{exact linearity of the LCG}. On the other hand, Harada et al.'s monotonous rhythm curves are  from ``almost linearity" of LCH, the continuum limit of LCG. 
Accordingly, shape analysis of various objects \cite{Harada1999, Kanaya2003en, Harada2001} and shape design \cite{Miura2010, Miura2014}  based on planar curves with ``almost linear LCG segments" were demonstrated.
Based on this background, we consider the LCG \eqref{LCG} of affine W-curves. 

Some curves in the affine W-curve were shown that their LCGs are asymptotically linear. 
The gradient of the LCG is computed from \eqref{LCG} as 
\begin{align}\label{LCGgradient}
    %\frac{dY}{dX}=1-\frac{\rho\rho_{ss}}{\rho_s}=1+\frac{\rho}{\rho_t^2}\placket{\frac{\rho_t s_{tt}}{s_t}-\rho_{tt}},\quad \rho =\frac{1}{\kappa^\mathrm{E}},
    \frac{dY}{dX}=1-\frac{\rho\rho_{ss}}{\rho_s}=1-\frac{\kappa^\mathrm{sim}_\theta}{(\kappa^\mathrm{sim})^2},  \quad \rho =\frac{1}{\kappa^\mathrm{E}},
\end{align}
which is asymptotically $2/3$ for the parabola and $1/2$ for the logarithmic curve \cite{Gobi2014}. 
The gradient of the LCG of the B\'ezier curve of degree $p\geq 3$ near the inflection points was shown to be $1/(2-p)$ \cite{Yoshida2010}. 

\begin{figure}[h]
\centering
    \includegraphics[height=0.4\columnwidth]{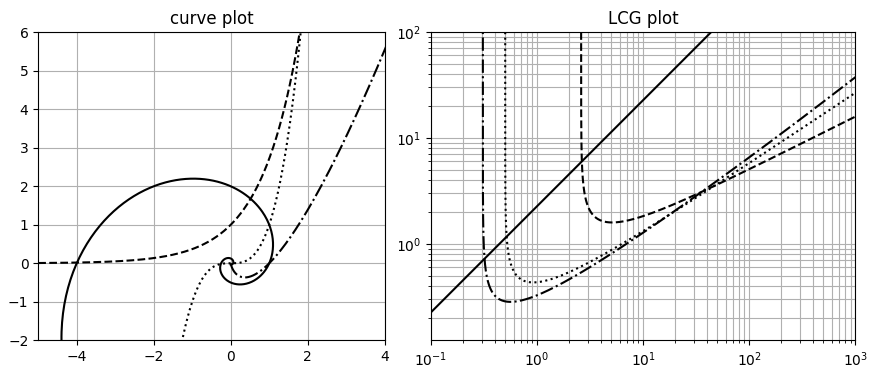}
       \caption{Some affine W-curves and their LCGs. Solid line:logarithmic spiral (slope $1$), dotted line: power function ($y=x^3$, asymptotic slope $3/5$), dashed line: exponential function (logarithmic function, LCG slope tends to $1/2$), dashdotted line: $x\log x$ (asymptotic slope $1$ as $x\to +\infty$). }
   \label{fig:fig2}
\end{figure}

\begin{figure}[h]
\centering
    \includegraphics[height=0.4\columnwidth]{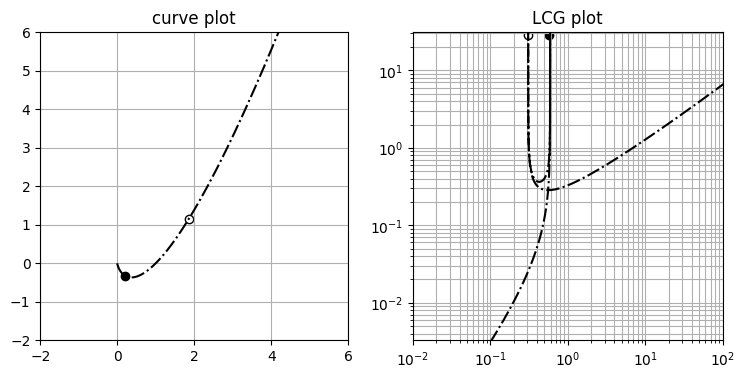}
       \caption{Three branches of the LCG of the graph of $y=x\log x$. The curve has no inflection point, and two critical (dotted) points of curvature, where the LCG branches. The asymptotic slope of both of the rightmost LCG-branch (as $x\to +\infty$) and the leftmost LCG-branch (as $x\to +0$) is $1$.}
   \label{fig:fig3}
\end{figure}

Figure \ref{fig:fig2} shows a single branch of the LCG of each of affine W-curves. 
As the logarithmic spiral is the LAC of slope $1$, its LCG is the line of slope $1$. 
For other affine W-curves (the parabola and the logarithmic curve, and the graph of $y=x\log x$), their LCGs are asymptotically linear. 
\begin{remark}
    We note that the LCG is subdivided into branches at the inflection point where $X=-\log\kappa^\mathrm{E}$ diverges, and at the critical points of curvature where $Y=\log|\kappa^\mathrm{E}/\kappa^\mathrm{E}_t|$ diverges. 
Figure \ref{fig:fig3} illustlates an example of the full branches of the LCGs and the corresponding curves. 

\end{remark}

\begin{figure}[h]
\centering
    \includegraphics[height=0.4\columnwidth]{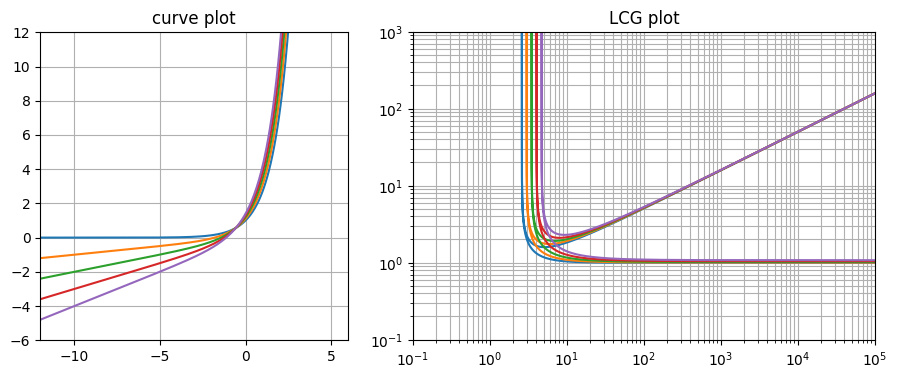}
       \caption{LCGs of one parameter family of affine deformations of graph of $y=e^x$. 
       %For each curve, the LCG branch neighbor to the inflection point (dot) is chosen for visibility.
       }
   \label{fig:fig6}
\end{figure}

\begin{figure}[h]
\centering
    \includegraphics[height=0.4\columnwidth]{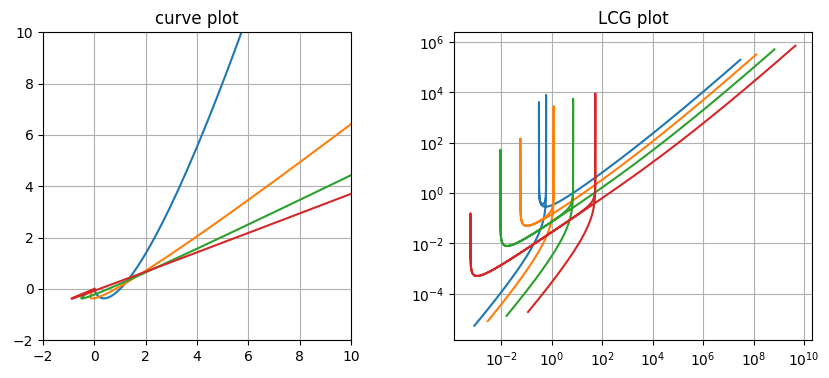}
       \caption{LCGs of one parameter family of affine deformations of graph of $y=x\log x$. 
       %For each curve, the LCG branch neighbor to the inflection point (dot) is chosen for visibility.
       }
   \label{fig:fig7}
\end{figure}

Basically, affine deformations do not preserve the LCGs of curves. 
Figure \ref{fig:fig6} shows 
affine deformations of the graph of exponential function, where the linearity of LCGs and their slopes are
preserved. 
With the exception of the logarithmic spiral, affine deformations of other W-curves
yield LCGs with asymptotically linear ends of the same slopes with originals  as shown in Figure \ref{fig:fig6} and \ref{fig:fig7}. 

\begin{figure}[h]
\centering
    \includegraphics[height=0.4\columnwidth]{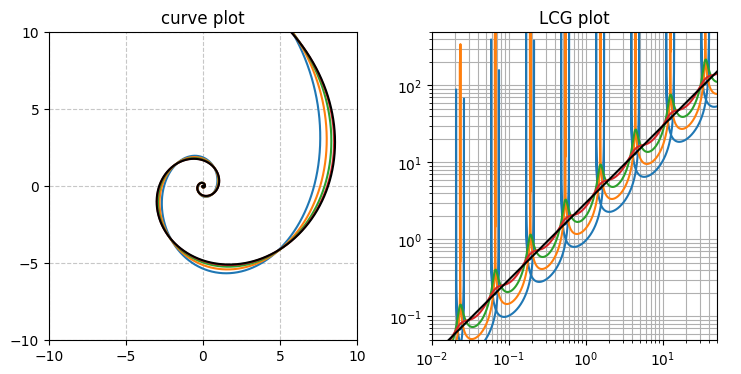}
       \caption{LCGs of one parameter family of affine deformations of  logarithmic spiral. }
   \label{fig:fig4}
\end{figure}

\begin{figure}[h]
\centering
    \includegraphics[height=0.37\columnwidth]{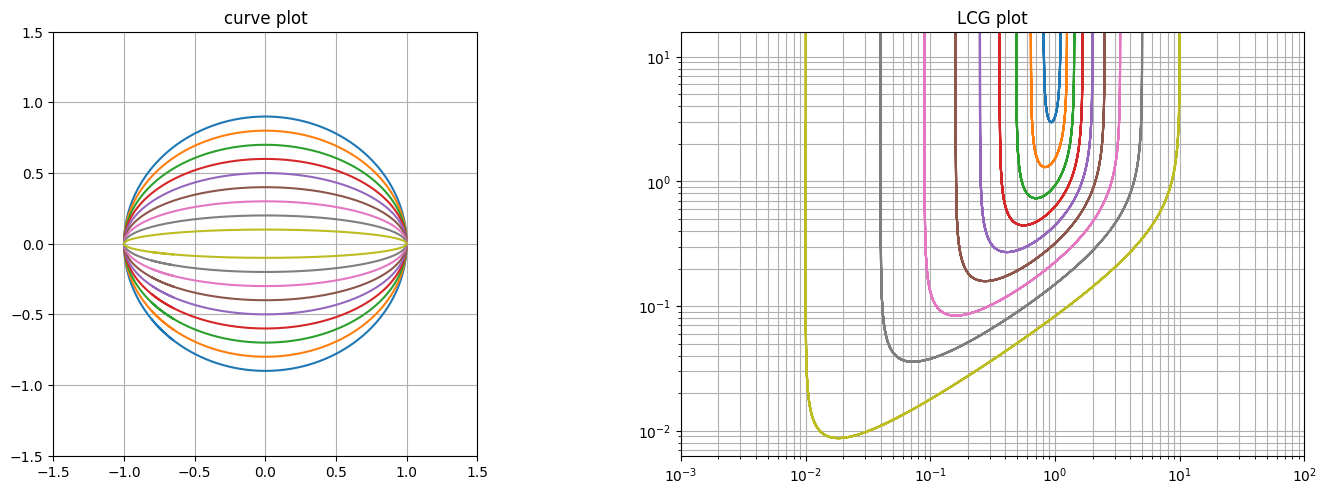}
       \caption{LCGs of one parameter family of affine deformations of ellipse. 
       As an ellipse approaches to the exact circle, LCG tends to a double vertical line, then, eventually degenerates to the infinity point. }
   \label{fig:ellipse}
\end{figure}

\begin{figure}[h]
\centering
    \includegraphics[height=0.4\columnwidth]{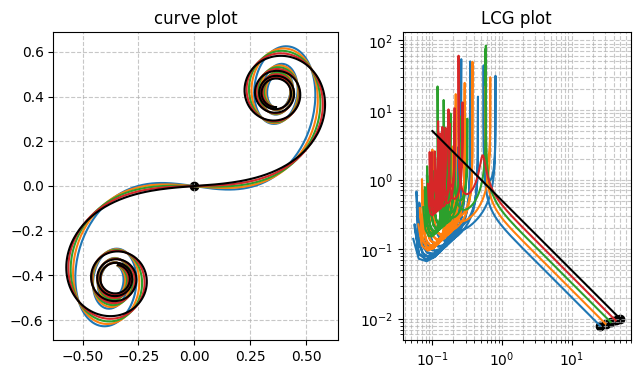}
       \caption{LCGs of one parameter family of affine deformations of the clothoid (Cornu spiral, LAC of slope $\alpha=2$). Near the center point, the LCGs share the linear tendency with common asymptotic slope $-1$. 
       Apart from the center, the LCGs yield many branches and their linear tendencies collapses, like the logarithmic spiral. }
   \label{fig:fig5}
\end{figure}

Figure \ref{fig:fig4} and Figure \ref{fig:ellipse} shows affine deformations of the logarithmic spiral and the ellipse, where the linearity of LCGs are collapsing as the curve deformed. 
We note about other LACs shown in Figure \ref{fig:fig5} as follows. 
Their LCGs also collapses but not repeatedly. Furthermore, they conain parts with the linear tendencies and slopes of the LCGs are preserved under affine deformations. 
These observations demonstrate that considering the affine W-curve and the LAC in a unified framework as a new family of ``aesthetic curves" may give us an important viewpoint.

%\note{placeholder text text text text text text text text text text text text text text text text text text text text text text text text text text text text text text text text text text text text text text text text text text text text text text text text text text text text text text text text text text text text text text text text text text text text}

\newpage
\section{Asymptotic Linearity of the LCGs for Affine W-curves}
In this section, we investigate the asymptotic behavior of the LCGs \eqref{LCG} for affine W-curves and establish the following theorem. 
Here we refer to the \textit{invariant directions of the ESA} of affine W-curves defined by the angles of the eigenvectors of Jacobian $DF(\varepsilon)\in \mathrm{GL}(2,\mathbb{R})$ in the ESA \eqref{sLACcurvest}. 
\begin{theorem}\label{thm:LCGslopeConst}
    %Let $\gamma$ be an affine W-curve other than a logarithmic spiral or an ellipse.     Then, for any turning angle $\theta_0$ corresponding to an asymptotic or inflection direction, the slope of the LCG asymptotically converges to a constant as $\theta \to \theta_0$. 
    Let $\gamma$ be a non-degenerate affine W-curve ($\kappa^\mathrm{E}\neq0$). %  except for affine deformations of the logarithmic spiral or the ellipse. 
    Let $\theta_0$ be an asymptotic turning angle corresponding to an invariant direction under the ESA. %or an inflection/singular point). 
    Then, the slope of the LCG asymptotically converges to a constant:
    \begin{align}\label{asymptoticslope}
        \lim_{\theta \to \theta_0} \frac{dY}{dX} = \begin{cases}
            1 + \placket{\underset{{\theta=\theta_0}}{\mathrm{Res}}\,\kappa^\mathrm{sim}(\theta) }^{-1} &\placket{{\mathrm{ord}_{\theta=\theta_0}}\,\kappa^\mathrm{sim}(\theta)=-1}
            \\[2mm]
            1 &\placket{{\mathrm{ord}_{\theta=\theta_0}}\,\kappa^\mathrm{sim}(\theta)=-2}, 
        \end{cases}
    \end{align}
    which covers all possibilities. %where no other orders occur. 
    Furthermore, the asymptotic behavior as $\theta\to\theta_0$ is invariant under the action of $G^\mathrm{A}$. 
\end{theorem}
As a conclusion, affine W-curves excluding some 
exceptional curves, namely, a logarithmic spiral, a circle, a line and their affine deformations, exhibit ``almost linearity'' on their LCGs near these asymptotic directions, asymptotically behaving as LACs. 
In particular, this property is invariant up to affine deformations. 
%Note that for the five representatives in Theorem~\ref{main}, $\theta_0$ as in Theorem~\ref{thm:LCGslopeConst} does not exist for logarithmic spiral and ellipse, and otherwise belongs to $\{0,\pm\pi/2,\pi\}$. 

\subsection{Asymptotic Expansion of the Similarity Curvature}\label{sec5.1}
We first note that the similarity curvature of a LAC of slope $\alpha$ \eqref{LACeq} is strictly constant when $\alpha = 1$, and is given by ${1}/\{{(1 - \alpha)(\theta-\theta_0)}\}$ when $\alpha \neq 1$ \cite{Inoguchi2018, Inoguchi2023}.
The following lemma claims that the similarity curvature of the affine W-curve asymptotically behaves like the LAC, up to affine transformations.

\begin{lemma}\label{lem:sim_curvature_expansion}
    Let $\gamma$ be a non-degenerate curve in the canonical representatives of affine W-curves in Theorem~\ref{main}.
    %Let $\gamma$ be an affine W-curve and $\theta_0$ be an invariant direction under the ESA. 
    %Up to affine transformations, $\gamma$ can be represented by a canonical curve whose similarity curvature $\kappa^\mathrm{sim}(\theta)$ admits a simple or double pole at $\theta_0$. 
    Then, the similarity curvature $\kappa^\mathrm{sim}(\theta)$ of $\gamma$ admits a simple or double pole at any invariant direction $\theta_0$ of the ESA of $\gamma$. 
    Specifically, it satisfies the asymptotic expansion
    \begin{align}\label{asymkappasim}
        \kappa^\mathrm{sim}(\theta) =\dfrac{l(\theta_0)}{(\theta - \theta_0)^2}+\dfrac{m(\theta_0)}{\theta - \theta_0} + \mathcal{O}(1) \quad (\theta \to \theta_0),
    \end{align}
    where $l(\theta_0),m(\theta_0)\in\mathbb{R}$ do not vanish simultaneously. 
\end{lemma}
%We note that in Lemma~\ref{lem:sim_curvature_expansion}, any $\theta_0$ are aligned with the coordinate axes so that any asymptotic turning angle $\theta_0$ belongs to $\{0, \pm\pi/2, \pi\}$.
%Hereafter we may assume $\theta_0 = 0$ without loss of generality up to the action of the dihedral group $D_4$ generated by a $\pi/2$-rotation and a reflection, which keep $(\theta,\kappa^\mathrm{sim}(\theta))$ up to signs and the LCG-slope invariant. 
%Throughout the following discussion, all asymptotic limits and expansions refer to $\theta \to 0$.
\begin{proof}
    We denote $y'=dy/dx$ and $y''=dy'/dx$. 
    For any canonical representative $\gamma$, the invariant directions are aligned with the coordinate axes so that any asymptotic turning angle $\theta_0$ belongs to $\{0, \pm\pi/2, \pi\}$.
    We may reduce to $\theta_0 = 0$ up to the action of the dihedral group $D_4$ generated by a $\pi/2$-rotation and a reflection, which keep $(\theta,\kappa^\mathrm{sim}(\theta))$ up to signs and the LCG-slope invariant. 
    %Throughout the following discussion, all asymptotic limits and expansions refer to $\theta \to 0$. 
    We discuss for each cases (a)-(e) according to Theorem \ref{main}. 
    \begin{enumerate}
        \item Let $\gamma$ be the graph of a power function $y = x^p$, $p \in \mathbb{R} \setminus \{0, 1\}$. 
        The turning angle satisfies $\tan \theta = y' = p x^{p-1}$, which yields $x = \left({\tan\theta}/{p}\right)^{\frac{1}{p-1}}$. 
        Differentiating $y'$ with respect to $x$ yields
        \begin{align}
            y'' = p(p-1) x^{p-2} = p(p-1) \left(\frac{\tan\theta}{p}\right)^{\frac{p-2}{p-1}},
        \end{align}
        from which the curvature radius is obtained as
        \begin{align}
            \rho(\theta) = \frac{\vlacket{1+(y')^2}^{3/2}}{|y''|} &= K (\cos\theta)^{-3} (\tan\theta)^{\frac{2-p}{p-1}} \nonumber
            \\&= K \theta^m \bigl(1 + \mathcal{O}(\theta^2)\bigr)\quad (\theta\to0),
        \end{align}
        where $m = {(2-p)}/{(p-1)}$ and $K > 0$ is a suitable constant. 
        Taking the logarithmic derivative, we obtain
        \begin{align}            \kappa^\mathrm{sim}(\theta) &= \frac{d\log\rho}{d\theta} %= \frac{m K \theta^{m-1} + \mathcal{O}(\theta^{m+1})}{K \theta^m \bigl(1 + \mathcal{O}(\theta^2)\bigr)} 
        = \frac{m}{\theta} + \mathcal{O}(\theta) \quad (\theta\to0).\label{ksimpow}      \end{align}
        
        For each $p$, the result \eqref{ksimpow} refers to the horizontal asymptotic end where $\theta_0=0$ (corresponding to $x \to 0$ for $p > 1$, and $x \to \infty$ for $p < 1$).  
        There is a vertical asymptotic end ($\theta_0=\pi/2$) remaining. 
        It is reduced to $\theta\to 0$ after applying the reflection interchanging the coordinate axes, which replaces $y = x^p$ with $y = x^{1/p}$, $\theta$ with $-\theta$, $\kappa^\mathrm{sim}$ with $-\kappa^\mathrm{sim}$, respectively. 
        Then, \eqref{ksimpow} results in the same presentation where $m$ is replaced with $({2 - {1}/{p}})({{1}/{p} - 1} )= ({2p - 1})/({1 - p})$. 
    %By the point or line symmetry of $\gamma$, 
    The above two evaluations essentially exhaust all asymptotic ends of the graph of $y=x^p$ for all $p \in \mathbb{R} \setminus \{0, 1\}$. 
    Thus, we have the conclusion with 
    \begin{align}
       ( l(0),m(0))=\placket{0,\frac{2-p}{p-1}},\quad (l(\pi/2),m(\pi/2))=\placket{0,\frac{2p - 1}{1 - p}}.
    \end{align}
        \item Let $\gamma$ be the graph of $y = \log x$. 
        Then, we have $\tan \theta = y'=1/x$, which implies 
            $x = 1/\tan\theta$. 
        Differentiating $y'$ with respect to $x$ gives $y''= -{1}/{x^2} = -\tan^2\theta$. 
        The curvature radius is given by
        \begin{align}
            \rho(\theta) &= \frac{\vlacket{1+(y')^2}^{3/2}}{|y''|} = \frac{1}{\sin^2\theta \cos\theta} = \theta^{-2}\bigl(1 + \mathcal{O}(\theta^2)\bigr)\quad (\theta\to+0),
        \end{align}
        from which we have 
        \begin{align}
            \kappa^\mathrm{sim}(\theta) %= \frac{d\log \rho}{d\theta} %= \frac{-2\theta^{-3} + \mathcal{O}(\theta^{-1})}{\theta^{-2}\bigl(1 + \mathcal{O}(\theta^2)\bigr)} 
            =(\log\rho(\theta))_\theta= -\frac{2}{\theta} + \mathcal{O}(\theta)\quad (\theta\to+0).
            \label{ksimlog}% \quad (\theta \to 0).
        \end{align}
        The result \eqref{ksimlog} refers to the horizontal asymptotic end ($\theta_0=0$, $x\to +\infty$).  
        The remaining end is vertical asymptotic ($\theta_0=\pi/2$, $x\to+0$), which is reduced to $\theta\to+0$ by considering $y=e^x$ instead. 
        Then, we have $y=y'=\tan\theta=y''=e^x$, which yields
        \begin{align}
            \rho(\theta)&= \frac{1}{\sin\theta \cos^2\theta} = \theta^{-1}\bigl(1 + \mathcal{O}(\theta)\bigr),
            \\
            \kappa^\mathrm{sim}(\theta)&=-\frac{1}{\theta}+\mathcal{O}(1). 
        \end{align}
        These two evaluations essentially exhaust all asymptotic ends of graph of $y=\log x$. 
    Thus, we have the conclusion with 
    \begin{align}
       ( l(0),m(0))=\placket{0,-2},\quad ( l(\pi/2),m(\pi/2))= \placket{0,-1}.
    \end{align}
        \item Let $\gamma$ be the graph of $y = x \log x$. 
        The original curve possesses two invariant asymptotic ends, the origin ($x \to +0$, $\theta_0 = -\pi/2$) and the infinite end ($x \to +\infty$, $\theta_0 = \pi/2$).
    Near the origin, applying the $\pi/2$-rotation in $D_4$ yields local representation $x = -y \log y$ with $y \to +0$, which is consistent with $\theta\to +0$ by $\cot \theta=1/y' = -(1+\log y)$ and $y=e^{-\cot\theta-1}$. 
        Multiplying $y'$ and the differential of $y'$ with respect to $y$ gives 
        \begin{align}
            y''=\frac{dy'}{dy}\frac{dy}{dx} = \frac{- 1/y}{(1+\log y)^2}\cdot\frac{1}{-(1+\log y)}={e^{\cot \theta+1}}{\tan^3\theta},
        \end{align}
        from which the curvature radius is given by
        \begin{align}
            \rho(\theta)&= \frac{\vlacket{1+(y')^2}^{3/2}}{|y''|} =\frac{e^{-\cot\theta-1}}{\sin^3\theta}.\label{rhoxlogx}
        \end{align}
        Taking the logarithmic derivative of \eqref{rhoxlogx} with respect to $\theta$, we have
        \begin{align}
            \kappa^\mathrm{sim}(\theta) =  \frac{1}{\sin^2\theta}-3\cot\theta =\frac{1}{\theta^2}-\frac{3}{\theta}+\mathcal{O}(1)\quad(\theta\to+0).\label{ksimxlogx}
        \end{align}
        The rest infinite end is reduced to $\theta\to +0$ by applying reflection in $D_4$ which replaces $x$ by $-x$, $\theta$ with $-\theta$, $\kappa^\mathrm{sim}$ with $-\kappa^\mathrm{sim}$, respectively. 
        Then, \eqref{ksimxlogx} results in 
        \begin{align}
            \kappa^\mathrm{sim}(\theta) =-\frac{1}{\theta^2}-\frac{3}{\theta}+\mathcal{O}(1)\quad(\theta\to+0).
        \end{align}
            Thus, we have the conclusion with 
    \begin{align}
       ( l(-\pi/2),m(-\pi/2))=\placket{1,-3},\quad ( l(\pi/2),m(\pi/2))=\placket{-1,-3}.
    \end{align}
        \item For a logarithmic spiral, there is no invariant direction $\theta_0$ under the ESA \eqref{ESAlogspiral}. 
        \item The parabola $y = x^2$ and the hyperbola $y = x^{-1}$ are already covered as graphs of power functions (a). 
        For a circle, there is no invariant direction $\theta_0$ under the ESA represented by rotations.  
    \end{enumerate}
    This completes the proof.
\end{proof}

%\subsection{Asymptotic Linearity and Slope of the LCG}
Based on the asymptotic behavior of the similarity curvature established in Lemma~\ref{lem:sim_curvature_expansion}, we now evaluate the asymptotic behavior of the LCG-slope. 
We may assume $\theta_0=0$ up to similarity transformations without loss of generality. 
Noticing %$\rho\dfrac{d}{ds}=\dfrac{d}{d\theta}$
$d/d\theta=\rho\cdot d/ds$, we have
\begin{align}
    \kappa^\mathrm{sim} &=(\log\kappa^\mathrm{E})_\theta=-(\log\rho)_\theta %= \frac{1}{\rho}\frac{d\rho}{d\theta} = \frac{d\rho}{ds} 
    = -\rho_s,
    \\
    \kappa^\mathrm{sim}_\theta&=\rho\kappa^\mathrm{sim}_s=-\rho\rho_{ss}.
\end{align}
%Then, the expression \eqref{LCGgradient} of LCG is rewritten by 
%\begin{align}\label{eq:dYdX_chain}
        %\frac{dY}{dX}         =1 +\frac{\rho\rho_{ss}}{(\rho_s)^2}        =1-\frac{\kappa^\mathrm{sim}_\theta}{(\kappa^\mathrm{sim})^2}=1-\frac{(\log|\kappa^\mathrm{sim}|)_\theta}{\kappa^\mathrm{sim}}. 
%    \end{align}
In the setting of Lemma~\ref{lem:sim_curvature_expansion}, substituting \eqref{asymkappasim} to 
%\eqref{eq:dYdX_chain} 
\eqref{LCGgradient} yields 
\begin{align}\label{eq:dYdX}
    \frac{dY}{dX}
    &=1-\frac{-\frac{2l}{\theta^3}-\frac{m}{\theta^2} + \mathcal{O}(1)}{\frac{l^2}{\theta^4}+\frac{2lm}{\theta^3} +\frac{m^2+lc}{\theta^2} + \mathcal{O}(\theta^{-1})}
    \\
    &\to\begin{cases}
        1+\dfrac{1}{m}&(l=0,m\neq0)\\
        1&(l\neq0)
    \end{cases} \quad (\theta\to0),
\end{align}
where $c\in\mathbb{R}$ is suitable constant. 
Then, we conclude as follows. 
\begin{proposition}\label{prop:ksimLCG}
    The asymptotic expansion \eqref{asymkappasim} of $\kappa^\mathrm{sim}$ implies the asymptotic linearity of LCG \eqref{asymptoticslope}. 
\end{proposition}
%we have Theorem \ref{thm:LCGslopeConst} for the five representatives (a)-(e) in Theorem \ref{main}. 
\subsection{Affine invariance of asymptotic LCG-slope}\label{sec5.2}
Finally, we discuss the invariance of the asymptotic LCG linearity under affine transformations. 
Let $F\in G^\mathrm{A}$ be an arbitrary affine transformation. 
Although quantities associated with Euclidean geometry $s,\kappa^\mathrm{E}$ and similarity geometry $\theta,\kappa^\mathrm{sim}$ are not preserved by $F$, the asymptotic behavior of the LCG remains invariant. 
Indeed, as is known in \textit{contact geometry} and \textit{jet theory} (see e.g.\,\cite{Olver}), non-singular affine transformations act linearly on tangent and contact spaces, thereby preserving the contact order and the characteristic scaling exponents of local asymptotic expansions.
Here, we present and prove the following  proposition without any additional setup. 
\begin{proposition}\label{prop:affineasymptotic}
    Let $\gamma$ be a regular curve, $F\in G^\mathrm{A}$ be an arbitrary affine transformation with Jacobian $DF\in \mathrm{GL}(2,\mathbb{R})$. 
    Denote the turning angle parameters of $\gamma,F\gamma$ by $\theta_\gamma,\theta_{F\gamma}$, and assume that the tangent direction at $\theta_\gamma \to 0$ is not mapped to the vertical asymptote under $DF$.
    Then, the turning angle $\theta_{F\gamma}$ and the Euclidean curvature $\kappa^\mathrm{E}_{F\gamma}$ satisfy
    \begin{align}
        \theta_{F\gamma} &= \theta_{F\gamma}^0 +c_1 \theta_\gamma + \mathcal{O}(\theta_\gamma^2)&& (\theta_\gamma \to 0), \label{thetaF}
        \\
        \kappa^\mathrm{E}_{F\gamma} &= c_2 \kappa^\mathrm{E}_\gamma \bigl(1 + \mathcal{O}(\theta_\gamma)\bigr) && (\theta_\gamma \to 0),
    \end{align}
    for some constants $\theta_{F\gamma}^0 \in \mathbb{R}$ and $c_1, c_2 \in \mathbb{R} \setminus \{0\}$ uniquely determined by $DF$.
\end{proposition}
\begin{proof}
    Throughout the following discussion, all asymptotic expansions refer to $\theta_\gamma \to 0$. 
    Let $M_Fz:={(c+dz)}/{(a+bz)}$ denote the M\"obius transformation \eqref{Mobslope} represented by the Jacobian $DF\in\mathrm{GL}(2,\mathbb{R})$. 
    By the assumption, we have $a\neq0$ and $\tan\theta_\gamma = \theta_\gamma + \mathcal{O}(\theta_\gamma^3)$. 
    Then, expanding $M_F\tan\theta_\gamma$ and its arctangent around $\theta_\gamma = 0$ yields 
    \begin{align}
        \tan \theta_{F\gamma} &=M_F\tan\theta_\gamma = \frac{c}{a} + \frac{ad - bc}{a^2} \theta_\gamma + \mathcal{O}(\theta_\gamma^2),
        \\
        \theta_{F\gamma}%&=\tan^{-1}\placket{\frac{c}{a} + \frac{ad - bc}{a^2} \theta_\gamma + \mathcal{O}(\theta_\gamma^2)}\\[2mm]
        &=\tan^{-1}(M_F\tan\theta_\gamma)=\tan ^{-1}\frac{c}{a}+\frac{ad-bc}{a^2+c^2}\theta_\gamma + \mathcal{O}(\theta_\gamma^2),
    \end{align}
    which is \eqref{thetaF}. 
    %where we have used $\tan^{-1} (k+l\theta_\gamma+\mathcal{O}(\theta_\gamma^2))=\tan^{-1}k+l\theta_\gamma/(1+k^2)+\mathcal{O}(\theta_\gamma^2)$. 
    \if0
    \begin{align}
        \tan \theta_{F\gamma} =M_F\tan\theta_\gamma = \frac{c}{a} + \frac{ad - bc}{a^2} \theta_\gamma + \mathcal{O}(\theta_\gamma^2).
    \end{align}
    Setting $\tan \theta_{F\gamma}^0 = c/a$, we have 
    \begin{align}
        \theta_{F\gamma} - \theta_{F\gamma}^0 &= \cos^2\theta_{F\gamma}^0¥ (\tan\theta_{F\gamma} - \tan\theta_{F\gamma}^0) + \mathcal{O}((\tan\theta_{F\gamma} - \tan\theta_{F\gamma}^0)^2)\nonumber\\
        &= c_1 \theta_\gamma + \mathcal{O}(\theta_\gamma^2), 
        \quad c_1 = \frac{\det(DF)}{a^2 + c^2} \neq 0.
    \end{align}
    \fi
    Next, we consider the transformation between the arc length parameters $s_\gamma$ of $\gamma$ and $s_{F\gamma}$ of $F\gamma$. 
    Noticing $(F\gamma(s_\gamma))_{s_\gamma}=DFe^{\imaginary\theta_\gamma}$, we have
    \begin{align}
        \frac{ds_{F\gamma}}{ds_\gamma} = |DF \cdot e^{\imaginary\theta_\gamma}| &= \sqrt{(a\cos\theta_\gamma + b\sin\theta_\gamma)^2 + (c\cos\theta_\gamma + d\sin\theta_\gamma)^2}
        \nonumber\\& = \lambda_0 \bigl(1 + \mathcal{O}(\theta_\gamma)\bigr), \quad \lambda_0 = \sqrt{a^2 + c^2} > 0.
    \end{align}
    Finally, from the chain rule we have
    \begin{align}
        \kappa^\mathrm{E}_{F\gamma} 
        = \frac{d\theta_{F\gamma}}{ds_{F\gamma}} 
        = 
        \frac{{d\theta_{F\gamma}}/{d\theta_\gamma}}{{ds_{F\gamma}}/{ds_\gamma}}\cdot\frac{d\theta_\gamma}{ds_\gamma}
        %{\frac{d\theta_{F\gamma}}{d\theta_\gamma} /\frac{d\theta_\gamma}{ds_\gamma}}{\frac{ds_{F\gamma}}{ds_\gamma}}
        = \frac{c_1 + \mathcal{O}(\theta_\gamma)}{\lambda_0 \bigl(1 + \mathcal{O}(\theta_\gamma)\bigr)} \kappa^\mathrm{E}_\gamma
        = c_2 \kappa^\mathrm{E}_\gamma \bigl(1 + \mathcal{O}(\theta_\gamma)\bigr),
    \end{align}
    where $c_2 = c_1 / \lambda_0 = {(ad-bc)}/{(a^2 + c^2)^{3/2}} \neq 0$. 
    This completes the proof.
\end{proof}
\begin{proof}[Proof of Theorem~\ref{thm:LCGslopeConst}]
We may assume $\theta_0=0$ up to similarity transformations $\gamma\mapsto R\gamma$, $R\in G^\mathrm{sim}$, where $\kappa^\mathrm{sim}_\gamma=\kappa^\mathrm{sim}_{R\gamma}$. 
Lemma~\ref{lem:sim_curvature_expansion} shows that any  affine W-curve is up to affine transformations represented by $\gamma$ whose similarity curvature has the principal part $l/\theta_\gamma^2+ m/\theta_\gamma$ around $\theta_\gamma=\theta_0$ for suitable $l,m\in\mathbb{R}$. 
It yields asymptotically linear LCG with slope expressible in terms of $l,m$ by Proposition~\ref{prop:ksimLCG}.  
By the expression \eqref{LCGgradient} and  Proposition~\ref{prop:affineasymptotic}, for any affine W-curve $F\gamma$, $F\in G^\mathrm{A}$, we have
\begin{align}
    \kappa^\mathrm{sim}_{F\gamma}=-\frac{d\log \kappa^\mathrm{E}_{F\gamma}}{d\theta_{F\gamma}}&=-\frac{\placket{\log\kappa^\mathrm{E}_\gamma+\log|c_2(F)|+\log(1+\mathcal{O}(\theta_\gamma))}_{\theta_\gamma}}{\placket{c_1(F)\theta_\gamma+\mathcal{O}(\theta_\gamma^2)}_{\theta_\gamma}}
    \\
    &=\frac{\kappa^\mathrm{sim}_\gamma+\mathcal{O}(1)}{c_1+\mathcal{O}(\theta_\gamma)}=\frac{\kappa^\mathrm{sim}_\gamma}{c_1}+\mathcal{O}(1)\quad (\theta_\gamma\to 0), \label{ksimF}
\end{align}
which ensures that $\kappa_\gamma^\mathrm{sim}(\theta_\gamma)$ having a double pole at $\theta_\gamma=0$ is equivalent to $\kappa_{F\gamma}^\mathrm{sim}(\theta_{F\gamma})$ having a double pole at $\theta_{F\gamma}=\theta_{F\gamma}^0$. 
If $\kappa^\mathrm{sim}_\gamma=m/\theta_\gamma+\mathcal{O}(1)\;(\theta_\gamma\to0)$, it follows from \eqref{thetaF} and \eqref{ksimF} that 
\begin{align}
    \kappa^\mathrm{sim}_{F\gamma}(\theta_{F\gamma})=\frac{m}{c_1\theta_\gamma}+\mathcal{O}(1)=\frac{m}{\theta_{F\gamma}-\theta^0_{F\gamma}}+\mathcal{O}(1)\quad (\theta_{F\gamma}\to \theta_{F\gamma}^0). 
\end{align}
Then, by Proposition~\ref{prop:ksimLCG}, the LCG-slope of $F\gamma$ coincides with that of $\gamma$. This completes the proof. 
\end{proof}
%We conclude Theorem~\ref{thm:LCGslopeConst} as follows. 
\begin{remark}
Here are some remarks on   Theorem~\ref{thm:LCGslopeConst}. 
\begin{itemize}
    \item In the above discussion of the proof, the assumption on $F$ in Proposition~\ref{prop:affineasymptotic} can be evaded by post-composition of similarity transformations $F\mapsto RF$ for $R\in G^\mathrm{sim}$, where $\kappa^\mathrm{sim}_{RF\gamma}=\kappa^\mathrm{sim}_{F\gamma}$. 
    %This replacement can be done independent of setting $\theta_0=0$
    \item The logarithmic spiral and ellipse have no invariant direction as mentioned in the proof of Lemma~\ref{lem:sim_curvature_expansion}. 
    The line is not non-degenerate since $\det(\gamma_s,\gamma_{ss})=\kappa^\mathrm{E}=0$. 
    These curves are exceptional curves for the claim. 
    \item Any asymptotic slope $\alpha\in\mathbb{R}$ can be attained by affine W-curves. 
    In fact, the LCG-slope of the graph of $y=x^p$, $p\in\mathbb{R}\setminus\{0,1\}$ is 
    \begin{align}
        \lim_{x\to 0}\frac{dY}{dX}&=1+\placket{\frac{2-p}{p-1}}^{-1}=\frac{1}{2-p},\\
        \lim_{x\to +\infty}\frac{dY}{dX}&=1+\placket{\frac{2p-1}{1-p}}^{-1}=\frac{p}{2p-1},
    \end{align}
    where the limit $p\to 0,\pm\infty $ is realized by the graph of $y=\log x$. 
    The LCG-slope $1$, attained by neither of them, is realized by the graph of $y=x\log x$. 
    \item On a regular point of the similarity curvature with $\kappa^\mathrm{sim}_\gamma(\theta)=\mathcal{O}(1)$ as $\theta\to\theta_1:=0$ (fixed up to $G^\mathrm{sim}$ without loss of generality), the LCG-slope is no longer affine invariant. 
    Indeed, substituting $\kappa^\mathrm{sim}_\gamma=a_0+a_1\theta_\gamma+\mathcal{O}(\theta_\gamma^2)$ to \eqref{eq:dYdX} yields 
    \begin{align}\label{eq:dYdXO(1)}
        \frac{dY}{dX}=1-\frac{a_1+\mathcal{O}(\theta_\gamma)}{(a_0+a_1\theta_\gamma+\mathcal{O}(\theta_\gamma^2))^2}\to 1-\frac{a_1}{a_0^2} \quad (\theta\to 0). 
    \end{align}
    The similar discussion as Proposition~\ref{prop:affineasymptotic} yields 
    \begin{align}
        \kappa^\mathrm{E}_{F\gamma} &= c_2 \kappa^\mathrm{E}_\gamma \bigl(1 +c_3\theta_\gamma+ \mathcal{O}(\theta_\gamma^2)\bigr) \quad (\theta_\gamma \to 0), 
    \end{align}
    for a suitable constant $c_3=c_3(F)\in\mathbb{R}\setminus\{0\}$. 
    Then, it follows from \eqref{ksimF} that
    \begin{align}
        \kappa^\mathrm{sim}_{F\gamma}&=-\frac{\placket{\log\kappa^\mathrm{E}_\gamma+\log|c_2(F)|+\log(1+c_3\theta_\gamma+\mathcal{O}(\theta_\gamma^2))}_{\theta_\gamma}}{\placket{c_1(F)\theta_\gamma+\mathcal{O}(\theta_\gamma^2)}_{\theta_\gamma}}\nonumber
    \\
    &=\frac{\kappa^\mathrm{sim}_\gamma-c_3}{c_1}+\mathcal{O}(\theta_\gamma^2)=\frac{a_0-c_3}{c_1}+\frac{a_1}{c_1}\theta_\gamma+\mathcal{O}(\theta_\gamma^2)\quad (\theta_\gamma\to 0), \\[2mm]
    &= \frac{a_0-c_3}{c_1}+\frac{a_1}{c_1^2}(\theta_{F\gamma}-\theta_{F\gamma}^0)+\mathcal{O}((\theta_{F\gamma}-\theta_{F\gamma}^0)^2)\quad (\theta_{F\gamma}\to\theta_{F\gamma}^0),
    \end{align}
    from which we have different asymptotic LCG-slope by \eqref{eq:dYdXO(1)}. 
\end{itemize}
\end{remark}

\section{Concluding discussions}
In this paper, we presented a family of curves characterized by the ESA in $G^\mathrm{A}$. 
It gives rise to a family of constant curvature curves in affine geometry, namely, the affine W-curve including the quadratic curve,  the logarithmic spiral, and the graphs of power functions, the logarithmic function, $y=x\log x$, and their affine equivalents.
On the other hand, the asymptotic linearity and their slopes of the LCGs of these curves are invariant under affine deformations, except for the logarithmic spiral and the circle. 
%Behavior of the LCG of the logarithmic spiral under affine deformation collapses like other LACs. 
%Based on the above observation, we propose a new family of ``aesthetic curves" consisting of 
As long as we regard the linear tendency of LCG as a characteristic of aesthetic shape element as Harada et al.\,suggested, it is natural to consider 

\begin{itemize}
    \item affine W-curves except for affine equivalents of a line, a circle, and a logarithmic spiral, 
    \item similarity equivalents of the logarithmic spiral and the exact circle, 
\end{itemize}
%On the other hand, LACs of slope $\alpha\neq 1$ admits the parts with the linear tendencies of the LCGs invariant under affine deformations, like other W-curves. 
as a new family of ``aesthetic curves", which shares the ESA in $G^\mathrm{A}$ as a common property. 

%Based on the ``almost linearity" of the LCG, we propose to consider the affine W-curve and the LAC in a unified framework as a new family of ``aesthetic curves". 
%Different from the HSA and the ESA, the MSA is a self-affinity along its similarity frame. 

On the other hand, the LAC is naturally described in the framework of integrable curve deformation theory in similarity geometry \cite{Inoguchi2018, Inoguchi2023} as an analog of Euler's elastica, which gives rise to an extension called the \textit{quasi-aesthetic curve (qAC)} \cite{InoguchiMiura2019, SatoSimizu2015, SatoShimizu2020}.
It includes the LAC, the quadratic curve, the elastic curve, and the graph of the logarithmic function.   
%\begin{itemize}
%    \item the LAC, 
%    \item the quadratic curve,
%    \item the elastic curve, and
%    \item the graph of the logarithmic function.  
%\end{itemize}
%Different from the HSA and the ESA, the MSA is a self-affinity along its similarity frame. 
It is interesting to note that the qAC overlaps with our family, despite the differences between the above two approaches. %, and the graph of the logarithmic function. 
The relationship between the affine W-curve and the qAC is yet to be discussed as a future work. 

Finally, %we remark on the MSA. 
since the similarity and affine transformation groups can be projected in subgroups of the {M\"obius} transformation group, the LAC and the affine W-curve together with self-affinities  may be unified in \textit{M\"obius geometry}. 
\if0
Further, substituting $t=t_0-\varepsilon$, equation \eqref{thetaSA} is rewritten by 
\begin{align}
\theta(t_0)
&=\frac{e^{-\frac{1-\alpha}{2}\varepsilon}\theta(t_0-\varepsilon)+e^{\frac{1-\alpha}{2}\varepsilon}\theta(t_0+\varepsilon)}{e^{-\frac{1-\alpha}{2}\varepsilon}+e^{\frac{1-\alpha}{2}\varepsilon}}. \label{inttheta}
\end{align}
Equation \eqref{inttheta} may have an application to engineering as an interpolation property arising from the MSA. 
\fi
We will work in these directions in the forthcoming publications.

\subsection*{Acknowledgments}
    We thank an anonymous referee for drawing the authors' attention to the Klein-Lie W-curve. 
    The authors would be grateful to Prof.\ Jun-ichi Inoguchi, Prof.\ Kenjiro T.\ Miura, Prof. Toshinobu Harada, and Prof.\ Yoshiki Jikumaru for their continuous encouragement. 
    This work was supported by JST CREST Grant Number JPMJCR1911, JSPS KAKENHI Grant Number 25K21661, 26K17041, and by Institute of Mathematics for Industry, Joint Usage/Research Center in Kyushu University.
    (FY2024 Workshop(I) “Evolving Design and Discrete Differential Geometry: towards Mathematics Aided Geometric Design” (2024a033))
\section*{Declarations}
\subsection*{Funding}
    This work was supported by JST CREST Grant Number JPMJCR1911, JSPS KAKENHI Grant Number 25K21661, 26K17041, and by Institute of Mathematics for Industry, Joint Usage/Research Center in Kyushu University.
    (FY2024 Workshop(I) “Evolving Design and Discrete Differential Geometry: towards Mathematics Aided Geometric Design” (2024a033))
\subsection*{Conflict of interest}
    The authors declare that there is no conflict of interest.
\subsection*{Ethics approval}
    Not applicable.
\subsection*{Consent to participate}
    Not applicable.
\subsection*{Consent for publication}
    The authors agree to publish this work.
\subsection*{CRediT authorship contribution statement}
    \textbf{Shun Kumagai:} Writing-Original draft preparation and Editing. 
    \textbf{Kenji Kajiwara:} Conceptualization, Writing-Reviewing and Editing.
\subsection*{Data availability}
    No data was used for the research described in the article. 
\subsection*{Code availability}
    Not applicable. 
 \bibliographystyle{elsarticle-num-names} 
\bibliography{reference.bib} 



\end{document}